\documentclass[preprint,12pt]{elsarticle}

\newcommand{\Su}[1]{{#1}^+ + {#1}^-}
\newcommand{\Di}[1]{{#1}^+ - {#1}^-}
\newcommand{\bx}{{\bf x}}
\newcommand{\by}{{\bf y}}
\newcommand{\bd}{{\bf d}}
\newcommand{\bn}{{\bf n}}
\newcommand{\I}{\mathrm{i}}
\newcommand{\tm}{\,\text{m}}
\newcommand{\um}{[\text{m}]}
\newcommand{\ud}{\text{d}}
\newcommand{\e}{\text{e}}

\usepackage{color}

\usepackage{graphicx}
\usepackage{amsmath,amssymb}
\usepackage{lineno}

\journal{Journal of Sound and Vibration}

\begin{document}
\begin{frontmatter}
  \title{Numerical simulation of sound propagation in and around ducts using thin boundary elements}
  \author{Wolfgang Kreuzer}
  \ead{wolfgang.kreuzer@oeaw.ac.at}


  \begin{abstract}
    Investigating the sound field in and around ducts is an important topic in acoustics, e.g. when simulating musical instruments or the human vocal tract. In this paper a method that is based on the boundary element method in 3D combined with a formulation for infinitely thin elements is presented. The boundary integral equations for these  elements are presented, and numerical experiments are used to illustrate the behavior of the thin elements. Using the example of a closed benchmark duct, boundary element solutions for thin elements and surface elements are compared with the analytic solution, and the accuracy of the boundary element method as function of element size  is investigated. As already shown for surface elements in the literature, an accumulation of the error along the duct can also be found for thin elements, but in contrast to surface elements this effect is not as big and a damping of the amplitude cannot be seen. In a second experiment, the impedance at the open end of a half open duct is compared with formulas for the radiation impedance of an unflanged tube, and a good agreement is shown. Finally, resonance frequencies of a tube open at both ends are calculated and compared with measured spectra. For sufficiently small element sizes frequencies for lower harmonics agree very well, for higher frequencies a difference of a few Hertz can be observed, which may be explained by the fact that the method does not consider dampening effects near the duct walls.  The numerical experiments also suggest, that for duct simulations the usual six to eight elements per wavelength rule is not enough for accurate results.
  \end{abstract}


\end{frontmatter}


\section{Introduction}
The simulation of the acoustics of tubes, ducts, and pipes has been of interest for some time now, and applications range from modeling musical instruments to the simulation of exhaust pipes or the human vocal tract, cf. \cite{CitLan11,Dubosetal99,FleRos91,Helieetal13,Pain05,Ruiz14,Wakita73,Webster19}.

The simplest and most common models 
assume the duct to have a constant (circular) cross section, and that its radius to wavelength ratio is small enough, so that plane wave propagation inside the duct can be assumed, cf. \cite{FleRos91}. Thus, the acoustic field can be approximated by a solution of the 1D Helmholtz equation or the 1D wave equation. In a segment with constant cross section this solution is given by a superposition of two plane waves traveling in opposite directions along the segment. At frequencies where the field of the two traveling waves add up to a standing wave resonances occur.  For specific boundary conditions at the duct ends (fully open or fully closed) the resonance frequencies for a duct with constant cross section can be calculated by a simple formula involving the length of the duct and the speed of sound (cf. \cite{FleRos91,Ruiz14}). For more realistic scenarios this simple formula can be modified by introducing correction factors for flanged and unflanged pipes based on the diameter of the duct, cf. \cite{FleRos91,Ruiz14,LevSch48,Silvaetal09}.

For non-constant cross sections the field can be calculated by dividing the duct into different segments of constant radius and by using the continuity of sound pressure and volume velocity across the segment borders (cf. \cite{Wakita73}). Alternatively, one can use the solution of Webster's horn equation for conical segments to derive the modes and thus the eigenfrequencies of the duct (cf. \cite{Webster19,Martin04}).

In case of open or half-open ducts, one problem with the simple 1D model is to correctly reproduce the radiation of the sound waves from the duct into the surroundings. To that end, finite element (FEM) or boundary element (BEM) methods can be used (cf. \cite{Marburg05}). For 2D and 3D simulations the FEM has become a popular tool (cf. \cite{Fussetal11,LefSca10,Marburgetal06,Matsuzakietal96,Vampolaetal08,Vampolaetal15}) because of its flexibility with respect to the geometry and the fact that the modes/resonances can be easily determined by solving a linear eigenvalue problem. However, for ducts with open ends the domain of interest is infinite and the FEM needs to be combined with infinite elements/absorbing boundary conditions/perfectly matched layers (PML) to correctly model the decay of the acoustic field towards infinity \cite{MarNol08}.

Compared to the FEM, the boundary element method has the advantage that only the surface of the scatterer needs to be discretized and that infinite domains needed for external scattering problems can be easily treated by this method. However, one drawback of the BEM is that the boundary integral equation (BIE) includes the Green's function for the Helmholtz equation, and especially for thin structures like plates or thin walled ducts the Green's function causes numerical problems, because the integrals involving both sides of the thin structure become nearly singular.

In this paper, we use a model for ducts, that is based on the BEM, but in contrast to most BEM models for tubes (for vocal tract simulation cf. \cite{Kagawaetal92, Motoki02}), the walls of the duct are modeled using thin elements described for example in \cite{Chenetal97,Krishnasamyetal90,Martinez91}. Using the formulation for thin elements has the advantage that some numerical problems with the almost singular integrals for very thin regular elements can be avoided.  Also, for sound hard or sound soft surfaces the formulation for the thin elements is very efficient.

To our knowledge, using BEM with thin elements has hardly been used in modeling the acoustics in and around thin walled ducts, and very little is known about the numerical behavior of these elements in that context. 
For regular surface elements it is, for example, known, that the accuracy of the BEM in connection with tubes is affected by a numerical damping effect (cf. \cite{BayMar18,Marburg16,Marburg18,Fahnline08}). However, as the BIE for thin elements is different, it is not clear, if and how such a numerical pollution effect influences the accuracy of a thin element formulation. Therefore, we will in a first numerical experiment test our implementation for thin elements with one of the benchmarks described in \cite{Hornikxetal15}, that has been extensively investigated for BEM using regular surface elements in \cite{BayMar18,Marburg16,Marburg18}. We briefly recall some of the results presented in \cite{BayMar18,Marburg16,Marburg18} and compare them with the results of an implementation using thin elements.

In a second experiment, we look at the radiation impedance for a half open duct that is modeled using thin elements and compare it with an approximation formula for unflanged tubes commonly used in literature \cite{FleRos91,Creasy16}. In a last example we calculate the resonance frequencies of a circular duct that is open at both ends (sound tube), compare the calculated resonance frequencies with measured spectra of the tube, and investigate the accuracy of the calculations with respect to the edge length of the elements used along the tube.

In all numerical experiments it is assumed that the sidewalls of the duct are acoustically hard and that the duct itself does not vibrate. We purely concentrate on the direct solution the Helmholtz equation to model sound propagation inside and outside the duct for a given set of frequencies. Although there is a considerable amount of research in connection with efficient methods to calculate eigenfrequencies and modes using the BEM (cf. \cite{Banerjeeetal88,Gaoetal13,Gaoetal19,NarBre83}), we leave that topic for future research.

The paper is organized as follows: First, we introduce the notation and some theory behind the 1D models and the BEM approach in Sections \ref{Sec:Notation} and \ref{Sec:Maths}. To keep the manuscript self contained we also present the formulation of the boundary integral equation used for the thin element BEM model. Numerical experiments with closed ducts, half open ducts, and a sound tube are made in Sections \ref{Sec:ExperimentI}, \ref{Sec:Impedance}, and \ref{Sec:Boomwhacker}. Finally, we conclude with a discussion of the results in Section \ref{Sec:Summary}. 

\section{Notation}\label{Sec:Notation}
In general, the sound field in the frequency domain at a point $\bx \in \mathbb{R}^3$ can be described by the velocity potential $\phi(\bx,k)$, where the wavenumber $k = \frac{2\pi f}{c}$ is given by the frequency $f$ and the speed of sound $c$ in the medium. To simplify notation, it is implicitly assumed that the velocity potential and all other related functions are always a function of the wavenumber, and $\phi(\bx)$ will be used instead of $\phi(\bx,k)$ most of the time.

From $\phi(\bx)$ the sound pressure $p(\bx)$ and the particle velocity $v(\bx,\bd)$ can then be derived using
\begin{align}\label{Equ:PresvsVeloPot}
  \nonumber p(\bx) = -\I \omega \rho \phi(\bx),\\
  v(\bx,\bd) = \frac{\partial \phi}{\partial \bd} = \nabla \phi \cdot \bd,
\end{align}
where $\bd$ is an arbitrary direction in space if $\bx$ lies in the free field, and $\bd = \bn$ if $\bx$ lies on the surface of an object. In this case, $\bn$ is the normal vector at $\bx$ pointing away from the object. $\nabla$ denotes the gradient $\nabla = (\frac{\partial}{\partial x},\frac{\partial}{\partial y},\frac{\partial}{\partial z})^T$, $\rho$ denotes the density of medium.

\begin{figure}[!h]
  \begin{center}
  \includegraphics[width=0.5\textwidth]{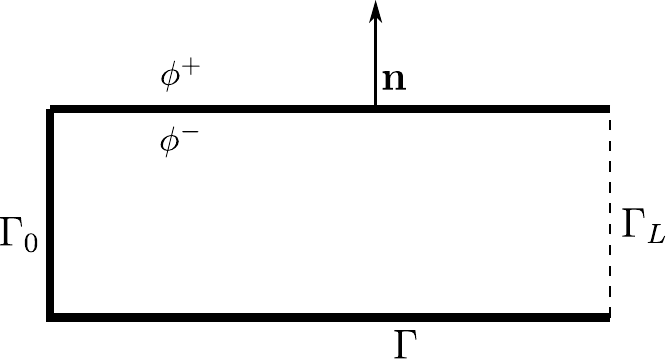}
  \caption{Example for the definition of the positive and negative sides for a half open duct. The sidewalls are denoted by $\Gamma$, the closure at $x = 0$ with $\Gamma_0$, at $x = L$ a (virtual) boundary $\Gamma_L$ is introduced.}\label{Fig:TubeScheme}
  \end{center}
\end{figure}

For the numerical experiments in Section~\ref{Sec:ExperimentI} admittance boundary conditions are used to model acoustic material properties at one duct ending. In this case, the particle velocity and the velocity potential are related using the scaled admittance $\alpha$:
$$
v(\bx) = \alpha(\bx) \phi(\bx).
$$
For (infinitely) thin middle face elements, the acoustic field on both sides of the element needs to be described. The positive side $\Gamma^+$ is defined by the direction of the normal vector (see for example Fig.~\ref{Fig:TubeScheme}). $\phi^+(\bx)$ denotes the velocity potential on this side, $\phi^-(\bx)$ is the velocity potential on the other side, the scaled admittances and particle velocities $\alpha^+, \alpha^-, v^+,$ and $v^-$ are defined similarly. Note, that for middle face elements the direction of the normal vector and thus the definition of positive and negative sides is arbitrary, in case of  ducts it seems natural to define the positive side as the outside and the negative side as the inside of the duct.

To simplify notation we follow \cite{Chenetal97} and introduce shortcuts for the sums and differences of the acoustic field on both sides:
\begin{align*}
  \bar{\alpha}(\bx) &= \Su{\alpha(\bx)}, \quad \phantom{\Delta}\bar{\phi}(\bx) = \Su{\phi(\bx)}, \quad \phantom{\Delta}\bar{v}(\bx) = \Su{v(\bx)},\\
  \Delta \alpha(\bx) &= \Di{\alpha(\bx)}, \quad \Delta \phi(\bx) = \Di{\phi(\bx)}, \quad \Delta v(\bx) = \Di{v(\bx)},
\end{align*}
for elements with admittance boundary condition we  set
\begin{align}
\nonumber   \bar{v}(\bx) = (\alpha^+ \phi^+ + \alpha^- \phi^-)(\bx) = \frac{1}{2}\left(\bar{\alpha} \bar{\phi} + \Delta \alpha  \Delta \phi\right)(\bx),\\
\label{Equ:Relations}  \Delta v(\bx) = (\alpha^+ \phi^+ - \alpha^- \phi^-)(\bx) = \frac{1}{2}\left(\Delta \alpha \bar{\phi} + \bar{\alpha}  \Delta \phi\right)(\bx).
\end{align}
\section{Mathematical background}\label{Sec:Maths}
\subsection{1D Duct Modeling}
Sound propagation inside a straight duct is an often used example for 1D acoustic wave propagation. If the width of the duct is small enough compared to the wavelength it can be assumed that the acoustic field inside the duct can be modeled using the 1D-Helmholtz equation, and the field is given by a superposition of two plane waves traveling from one end of the duct to the other in opposite directions. In this case, the solution (velocity potential $\phi$, pressure $p$, or particle velocity $v$) in the frequency domain can be written as
\begin{align}
  \nonumber \phi(x,k) &= A \e^{\I k x} + B \e^{-\I k x},\\
  \nonumber p(x,k) &= -\I \omega \rho \left(A \e^{\I k x} + B \e^{-\I k x}\right),\\
  \label{Equ:Helmholtz1D} v(x,k) &= \I k \left(A \e^{\I k x} - B \e^{-\I k x}\right),
\end{align}
where the harmonic time dependence is given by $\e^{-\I \omega t}$, thus the velocity potential for a harmonic tone with frequency $\omega$ in the time domain is given by $\check{\phi}(x,t) = \phi(x,k)\e^{-\I \omega t}$. The unknowns $A$, $B$, and possible resonance frequencies are determined by the boundary conditions at both ends of the duct. For the numerical experiments in Section \ref{Sec:ExperimentI}, for example, it is assumed, that both ends of the duct are closed. At $x = 0$  a Neumann boundary condition $v(0,k) = v_0$\,ms$^{-1}$ is used to define a sound source, at the other end at $x = L$ an impedance boundary condition $Z(L,k)  = \frac{p(L,k)}{v(L,k)}= \rho c$ is used. To fulfill these boundary conditions, $A$ needs to be zero, and the analytic solution $\phi(x,k)$ is given by a single plane wave traveling along the duct. 

\subsection{BEM and thin elements}
To model acoustic wave propagation in 3D in the frequency domain, one usually solves the (free field) Helmholtz equation
\begin{align}
\label{Equ:Helmholtz}  \nabla^2 \phi(\bx) + k^2 \phi(\bx) = \phi_{\text{inc}}(\bx),  \quad  \bx \in \Omega,\\
\label{Equ:Helmholtzbc}  v(\bx) = v_0(\bx) + \alpha(\bx) \phi(\bx), \quad \bx \in \Gamma,\\
\label{Equ:Sommerfeld}  \lim\limits_{r\rightarrow \infty} r  \left( \frac{\partial }{\partial r} - \I k\right) \phi(\bx) = 0,
\end{align}
where $r = ||\bx||$. Eq.~(\ref{Equ:Helmholtzbc}) denotes Robin boundary conditions on the surface $\Gamma$ of the scatterer, where $v_0$ is a possible known given boundary velocity and Eq.~(\ref{Equ:Sommerfeld}) represents the Sommerfeld radiation condition, that guarantees the decay of the acoustic field at infinity. $\Omega$ is either the domain inside or outside the scatterer, respectively. $\phi_{\text{inc}}(\bx)$ denotes the sound field of an external sound source, e.g., a point source or a plane wave.

For the BEM  a weak formulation of Eq.~(\ref{Equ:Helmholtz}) is combined with the Green's function for the free field Helmholtz equation 
$$
G(\bx,\by) = \frac{\e^{\I k ||\bx-\by||}}{4\pi ||\bx-\by||}
$$
as test function. For points $\bx$ on a smooth part of $\Gamma$ the boundary integral equation (BIE) is then given as 
\begin{equation}
  \frac{1}{2}\phi(\bx) - \tau\int\limits_\Gamma H(\bx,\by) \phi(\by) \ud\by + \tau\int\limits_{\Gamma} G(\bx,\by) v(\by) \ud\by = \int\limits_{\Omega} G(\bx,\by) \phi_{\text{inc}}(\by) \ud\by, 
\end{equation}
where $H(\bx,\by) = \frac{\partial G}{\partial \bn_y}(\bx,\by)$. $\bn_y$ is the normal vector pointing to the outside of the scatterer at a point $\by$. The factor $\tau$ is used to distinguish between an exterior problem ($\tau = 1$) where the field outside the scatterer is of interest, and interior problems ($\tau = -1$).

Using approaches like collocation, the Galerkin, or the Nystr\"om method \cite{Nystroem30,SauSch11} the BIE is transformed into a linear system of equations that needs to be solved for an approximation of the unknown velocity potential $\phi(\bx)$, depending on the method used. For the numerical experiments in this paper, we use the collocation method with constant elements. For that, the surface of the scatterer is discretized using triangular or quadrilateral elements  $\Gamma_j$. On each of these little elements the velocity potential is assumed to be constant with (yet) unknown potential $\phi_j$. The BIE is then evaluated at a given set of points (= collocation points), which in our case are given by the midpoints $\bx_j$ of each element.  This approach has the advantage that
\begin{itemize}
\item the application of the collocation method is relatively easy and efficient,
\item as the collocation point is positioned in the middle of the element, the velocity potential never needs to be evaluated at the edges of the elements, where singularities may occur, especially  at free edges,
\item as the discretization uses non-continuous elements, a refinement of the mesh in selected parts can be easily done without having to adapt the mesh in other parts.
\end{itemize}

In \cite{Chenetal97,KreWeb19} the equations for the thin elements were derived by looking at one ``standard surface'' BEM-element modeling a part of the boundary with fixed thickness $s>0$, by deriving the BIEs for this element and then by taking the limit of these equations for  $s\rightarrow 0$.

To keep the manuscript self contained, we will in the following derive the BIEs for thin elements, but use an alternative approach combining the formulations for exterior and interior problem, which, in our opinion, opens new views on the behavior of the numerical solution for ducts. We use the example of a half open duct of length $L$ which is open at $x = L$, where a  virtual boundary $\Gamma_L$ is introduced. The sidewalls are denoted by $\Gamma$, all elements on $\Gamma$ are assumed to be sound hard, thus, $v(\bx) = 0, \bx \in \Gamma$. At $x = 0$ the duct is closed by $\Gamma_0$, where a velocity boundary condition $v(\bx) = v_0, \bx \in \Gamma_0$ is assumed (cf. Fig.~\ref{Fig:TubeScheme}). To keep the notation simple, no external sound source is given. 

With the virtual closure $\Gamma_L$ the half open duct can be interpreted as the coupling of an internal problem and an external problem. For $\bx$ on a smooth part of the boundary $\tilde{\Gamma} = \Gamma \cup \Gamma_0$ the formulation for the \emph{interior} problem is given by 
\begin{align}
  \nonumber\frac{\phi^{-}(\bx)}{2} &+ \int\limits_{\tilde{\Gamma}} H(\bx,\by) \phi^{-}(\by)\ud\by - \int\limits_{\tilde{\Gamma}} G(\bx,\by) v^{-}(\by)\ud\by\\
    \label{Equ:Interior}                   &+ \int\limits_{\Gamma_L} H(\bx,\by) \phi^{-}(\by)\ud\by - \int\limits_{\Gamma_L} G(\bx,\by) v^{-}(\by)\ud\by = 0.
\end{align}
The ``$.^-$'' superscript  over the velocity potential $\phi(\bx)$ and the particle velocity $v(\bx)$ indicates that we are looking at the ``interior'' part of the duct. 

The BIE for the \emph{exterior} part is given as
\begin{align}
  \nonumber\frac{\phi^{+}(\bx)}{2} &- \int\limits_{\tilde{\Gamma}} H(\bx,\by) \phi^{+}(\by)\ud\by + \int\limits_{\tilde{\Gamma}} G(\bx,\by) v^{+}(\by)\ud\by \\
    \label{Equ:Exterior}        &- \int\limits_{\Gamma_L} H(\bx,\by) \phi^{+}(\by)\ud\by + \int\limits_{\Gamma_L} G(\bx,\by) v^{+}(\by)\ud\by = 0.
\end{align}
Combining Eqs.~(\ref{Equ:Interior}) and (\ref{Equ:Exterior}) and using the fact that on $\Gamma_L$ the sound field is the same on both sides ($\phi^+(\bx) = \phi^{-}(\bx), v^{+}(\bx) = v^{-}(\bx), \bx \in \Gamma_L$) yields
\begin{align}
\nonumber  \frac{  \phi^{+}(\bx) + \phi^{-}(\bx) }{2} &- \int\limits_{\tilde{\Gamma}}  H(\bx,\by) \left( \phi^+(\by) \!-\! \phi^-(\by) \right)\ud\by \;+\\
\label{Equ:Thin1}  &+\int\limits_{\tilde{\Gamma}} G(\bx,\by) \left( v^+(\by) \!-\! v^-(\by) \right)\ud\by = 0, 
\end{align}
or with the notation introduced in Section~\ref{Sec:Notation}
\begin{equation}\label{Equ:Midfirst}
  \frac{\bar{\phi}(\bx)}{2} - \int\limits_{\tilde{\Gamma}} H(\bx,\by) \Delta \phi(\by) \ud\by+ \int\limits_{\tilde{\Gamma}} G(\bx,\by) \Delta v(\by) \ud\by = 0.
\end{equation}
Using collocation with constant elements, the above equation can be transformed into a linear system of equations that has, in general,  about twice as many unknowns as equations unless pressure boundary conditions are used everywhere. In that special case $\Delta \phi$ and $\bar{\phi}$ are known, leaving $\Delta v$ as the only unknown. 
In all other cases an additional system of equations is needed. These additional equations are created by taking the derivative of Eq.~(\ref{Equ:Midfirst}) with respect to the normal vector $\bn_x$ at the point $\bx$ to get
\begin{equation}\label{Equ:Midsecond}
  \frac{\bar{v}(\bx)}{2} - \int\limits_{\tilde{\Gamma}} E(\bx,\by) \Delta \phi(\by) \ud\by+ \int\limits_{\tilde{\Gamma}} H'(\bx,\by) \Delta v(\by) \ud\by = 0,
\end{equation}
where $H'(\bx,\by) = \frac{\partial G}{\partial \bn_x}(\bx,\by)$ and $E(\bx,\by) = \frac{\partial^2 G}{\partial \bn_x \partial \bn_y} (\bx,\by)$. 

Let us now assume, that we can split $\tilde{\Gamma} = \Gamma \cup \Gamma_a \cup \Gamma_0$ into  a part $\Gamma$ with acoustically hard boundary conditions ($v^{\pm}(\bx) = 0, \bx \in \Gamma$), a part $\Gamma_a$ with admittance boundary conditions ($v^{\pm}(\bx) = \alpha^{\pm}(\bx) \phi^{\pm}(\bx), \bx \in \Gamma_a$) and a part $\Gamma_0$  with  a velocity boundary condition $v^{\pm}(\bx) = v_0^{\pm}, \bx \in \Gamma_0$.

In that setting the BIE can be formulated using 3 sets of equations. For the points with admittance boundary conditions ($\bx \in \Gamma_a)$, we get
\begin{align}
  \frac{\bar{\phi}(\bx)}{2} - \int\limits_{\tilde{\Gamma}} H(\bx,\by) \Delta \phi(\by) \ud\by + \int\limits_{\Gamma_a} G(\bx,\by) \Delta v(\by) \ud\by &= -\int\limits_{\Gamma_0} G(\bx,\by) \Delta {v_0} \ud\by, \\
  \frac{\bar{v}(\bx)}{2} - \int\limits_{\tilde{\Gamma}} E(\bx,\by) \Delta \phi(\by)\ud\by + \int\limits_{\Gamma_a} H'(\bx,\by) \Delta v(\by) \ud\by &= -\int\limits_{\Gamma_0} H'(\bx,\by) \Delta {v_0}\ud\by
\end{align}
or by using Eq.~(\ref{Equ:Relations})
  \begin{multline}\label{Equ:BIEmMidII} 
    \frac{\bar{\phi}(\bx)}{2} - \int\limits_{\Gamma} H(\bx,\by) \Delta {\phi}(\by)\ud\by + \int\limits_{\Gamma_a} \frac{G(\bx,\by)}{2}
    \left(
        \Delta \alpha \bar{\phi} + \bar{\alpha}\Delta \phi
    \right)(\by)
     \ud\by = \\
     - \int\limits_{\Gamma_0} G(\bx,\by) \Delta {v_0} \ud\by,\quad \bx \in \Gamma_a,
  \end{multline}
and
  \begin{multline}\label{Equ:BIEmidIII}
    \frac{(\bar{\alpha} \bar{\phi} + \Delta \alpha \Delta \phi)(\bx)}{4} - \int\limits_{\Gamma} E(\bx,\by) \Delta {\phi}(\by)\ud\by  + \int\limits_{\Gamma_a} \frac{H'(\bx,\by)}{2}
    \left(
    \Delta \alpha \bar{\phi} + \bar{\alpha}\Delta \phi
    \right)(\by)
    \ud\by \\= - \int\limits_{\Gamma_0} H'(\bx,\by) \Delta {v_0} \ud\by ,\quad \bx \in \Gamma_a.
  \end{multline}
For the rest of the points it is sufficient to use the equation
\begin{multline}\label{Equ:BIEmidI}
  \int\limits_{\Gamma} E(\bx,\by) \Delta \phi(\by) \ud\by - \int\limits_{\Gamma_a} \frac{H'(\bx,\by)}{2} \left(
  \bar{\alpha} \Delta {\phi} + \Delta {\alpha} \bar{\phi}
  \right) (\by) \ud\by = \\
 \int\limits_{\Gamma_0} H'(\bx,\by) \Delta {v_0} \ud\by + \frac{\bar{v_0}}{2}\delta_{\Gamma_0}(\bx),
\end{multline}
where $\delta_{\Gamma_0}(\bx) = 1$ if $\bx \in \Gamma_0$ and $\delta_{\Gamma_0}(\bx) = 0$ otherwise. Thus, the system consists of $n + n_0 + 2\cdot n_a$ unknowns, ($\Delta \phi$ for each of the $n$ elements on $\Gamma$, $\Delta \phi$ for each of the $n_0$ elements on $\Gamma_0$, and $\Delta \phi$ and $\bar{\phi}$ for each of the $n_a$ elements on $\Gamma_a$) and the same number of linear equations.

Once the above system of equations is solved, $\bar{\phi}$ for $\bx \in \Gamma \cup \Gamma_0$ can be simply determined by inserting the solution of the above system into
\begin{multline}
  \frac{ \bar{\phi}(\bx) }{2} - \int\limits_{\tilde{\Gamma}}H(\bx,\by) \Delta \phi(\by) \ud\by + \int\limits_{\Gamma_a} G(\bx,\by)
  \left(
     \bar{\alpha} \Delta {\phi} + \Delta {\alpha} \bar{\phi}
  \right) (\by) \ud\by \\
  = -\int\limits_{\Gamma_0} G(\bx,\by) v_0 \ud\by.
\end{multline}
\section{Numerical experiment: Closed duct}\label{Sec:ExperimentI}
The  first numerical experiment assumes a duct with a square cross section with width $w = 0.2$\,m and a length $L = 3.4$\,m. Both ends are closed, at the closure $\Gamma_0$  at $x = 0$ a velocity boundary condition with $v_0 = (1,0,0)^T$\,ms$^{-1}$ is applied, and at the other end at $x = L$ an impedance of $Z = \frac{p}{v} = \rho c$ is assumed, where $\rho$ and $c$ are the density of the air and the speed of sound, respectively. The sidewalls $\Gamma$ of the duct are assumed to be sound hard $v(\bx) = 0$\,ms$^{-1}$, $\bx \in \Gamma$. Using  Eq.~(\ref{Equ:Helmholtz1D}) the analytic solution of the Helmholtz equation for this duct is given by $\phi(x,k) = \frac{1}{\I k} \e^{-\I kx}$. In all experiments in this and the following sections, the density $\rho$ will be set to $\rho = 1.3$\,kg\,m$^{-3}$, the speed of sound will be assumed to be $c = 340$\,m\,s$^{-1}$. 
This example is one of the benchmark problems in \cite{Hornikxetal15} and BEM calculations with regular surface elements have been thoroughly tested in \cite{BayMar18,Marburg18}.

In the following we a) partly reconstruct the results from the numerical experiments in \cite{Marburg18} for regular surface elements, and b) compare these results to numerical experiments using thin elements.
We are aware that the nature of the benchmark as an interior problem is not exactly suited for thin elements, where the sound field on both sides of the element is of interest, but this benchmark has a simple analytic solution and the calculated field $\phi^+$ on the outside walls of the duct, which should be zero in theory, provides an additional error indicator for the thin elements. 

\begin{figure}[!h]
  \begin{center}
    \includegraphics[width=0.59\textwidth]{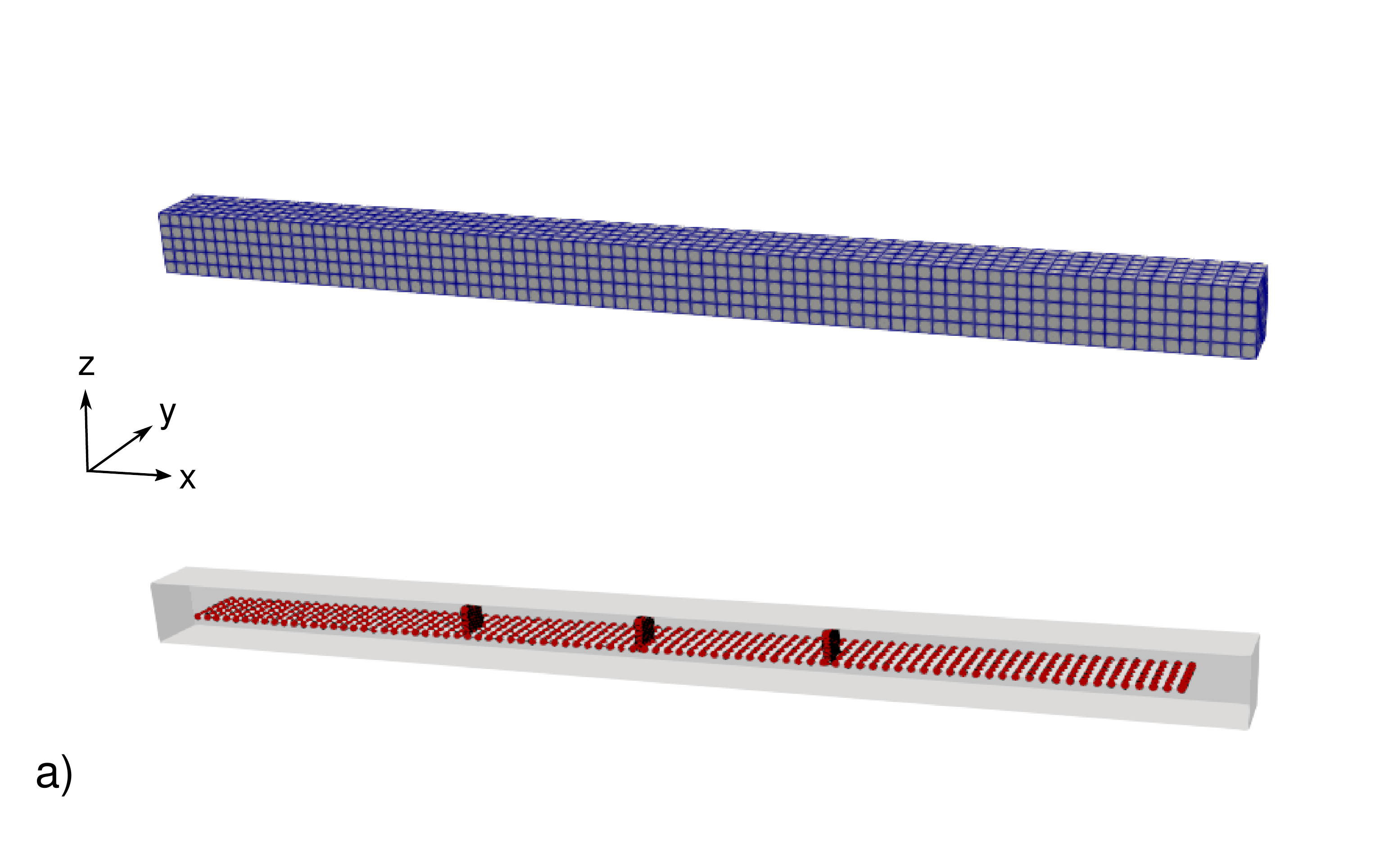}
    \includegraphics[width=0.39\textwidth]{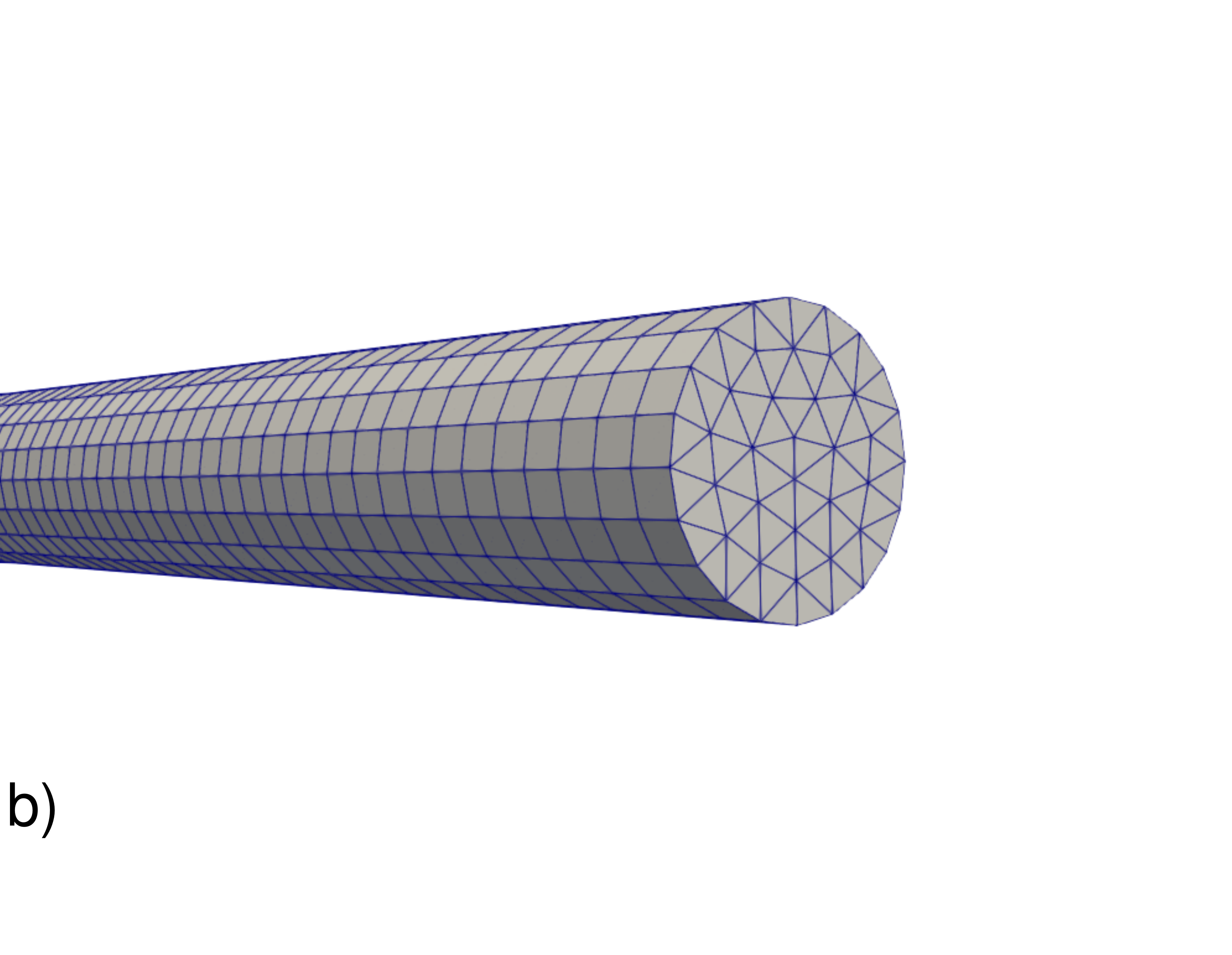}
  \end{center}
  \caption{a) Standard  mesh for a square duct and the positions of the evaluation points inside a square duct duct. b) Discretization of the end of a circular duct.}\label{Fig:Mesh}
\end{figure}

The default size of the quadrilateral BE-elements used to discretize the sidewalls was $h_1 \times h_2 = 0.04$\,m$\times 0.04\,$m, which is equivalent to a discretization using 20 elements around the duct and 85 elements along its length. The closures of the duct were discretized with elements of size $h_0 \times h_0 = h_L \times h_L = 0.04\tm \times 0.04\tm$. For the sound field inside the duct, Green's representation formula
\begin{equation}\label{Equ:GreensRepresentation}
  \phi(\bx_i) = \int\limits_{\Gamma_{\text{all}}} H(\bx_i,\by) \phi(\by)\ud\by - \int\limits_{\Gamma_{\text{all}}} G(\bx_i, \by) v(\by) \ud\by
\end{equation}
was used to evaluate the sound field at points $\bx_i = (x_i,y_i,z_i)$ inside the duct positioned at:
\begin{itemize}
  \item a $7 \times 86$ grid in a plane along the middle of the duct where $x_i$ lies between 0.17\,m and 3.23\,m, $y_i$ lies between $-0.09$\,m and 0.09\,m, and $z_i$ is fixed to  0\,m,
  \item a $16 \times 16$ grid in planes normal to the $x$-axis at $x_i \in \{L/3,L/2,2L/3\}$ where $y_i$ lies between 0\,m and 0.08\,m and $z_i$ lies between 0\,m and 0.08\,m.
\end{itemize}
In Eq.~(\ref{Equ:GreensRepresentation}) $\Gamma_{\text{all}} = \Gamma \cup \Gamma_0 \cup \Gamma_L$ is the union of the sidewalls and the closures of the duct (see Fig.~\ref{Fig:TubeScheme} for the definitions). The mesh of the duct and the position of the evaluation points inside the duct are depicted in Fig.~\ref{Fig:Mesh}a. In the next subsection BEM calculations for a square duct are compared with calculations for a circular duct with the same cross sectional area. The mesh for this duct is shown in Fig~\ref{Fig:Mesh}b.

\subsection{Duct with regular surface elements}
\subsubsection{Field along the duct}\label{Sec:Surface}
\begin{figure}[!h]
  \begin{center}
    \includegraphics[width=0.95\textwidth]{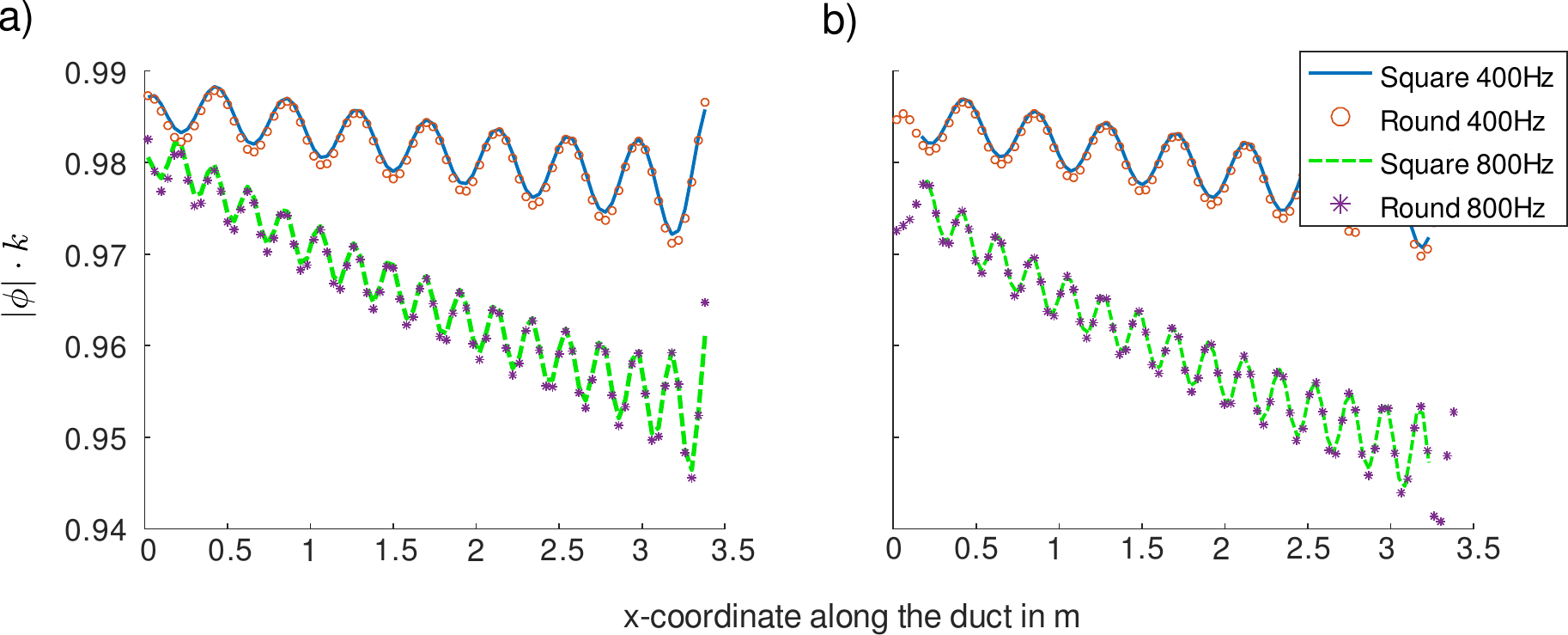}
   \end{center}
   \caption{a) Absolute values of the scaled velocity potential along a line at the side of the duct and b) the scaled velocity potential inside the duct along its center. Lines depict the velocity potential $\phi$ for the square duct, the markers are used for the circular duct. The continuous line shows the velocity potential for $f =  400$\,Hz, the dashed line the values for $f = 800$\,Hz. All calculations have been done using surface elements.}\label{Fig:RoundvsSquare}
\end{figure}

As long as the cross sectional area stays the same, the shape of the cross section does not have much impact on the calculated velocity potential except maybe at the corners. In Fig.~\ref{Fig:RoundvsSquare}a the (scaled) velocity potentials at collocation points along a line on a sidewall of  two ducts with circular and square cross sections are compared. The square duct had a width of $w = 0.2\tm$; in order to have the same cross-sectional area the radius of the circular tube was set to $r = \frac{w}{\sqrt{\pi}}$. 
For a better comparison of the curves for the frequencies $f = 400\,$Hz and $f = 800$\,Hz the velocity potentials were scaled with the wavenumber $k$.

For the square duct the velocity potential on the duct walls was evaluated at collocation points on a line along the middle of one side at points $\bx_i = (x_i, y_i, z_i)$ with $0.0\tm < x_i < 3.4\tm$, $y_i = -0.1\tm$, and $z_i = 0.0\tm$. As the field on the circular duct is rotational symmetric, any line along the duct can be used as reference line, for technical reasons $y_i \approx 0.11\tm$ and $z_i \approx 0.017\tm$, the $x$-coordinates of the points were the same as for the square duct.

The square duct was discretized using the default rectangular elements of size $h_1 \times h_2 = 0.04\tm \times 0.04\tm$, for the circular duct the procedure to create the mesh was a bit more involved. First, the discs closing the right and left sides of the duct were discretized with triangular elements using the matlab/octave toolbox \texttt{distmesh} \cite{PerStr04}. The mean edge length of the triangles was chosen such that 20 elements are used along the circumference of the disc,  resulting in a mean edge length of about $0.035\tm$. The sidewalls of the duct were discretized using rectangular elements with a fixed width $h_1 = 0.04\tm$ and slightly varying height $h_2 \approx 0.035\tm$ to match the discretization of the discs at the closures, see also Fig.~\ref{Fig:Mesh}b, where the mesh of the circular duct is displayed. The main characteristics of the different meshes used in this subsection are summarized in Tab.~\ref{Tab:SummarySurface}.

As in \cite{Marburg18,Fahnline08}, a numerical damping effect for the acoustic field can be observed in Fig.~\ref{Fig:RoundvsSquare}, that depicts the scaled absolute value of the velocity potential $\phi$ on the surface as well as the inside the duct, where the velocity potential was evaluated at a line along the middle of the duct. The continuous and dashed lines depict $|\phi|\cdot k$ for the square duct at 400\,Hz and 800\,Hz, respectively. The '$\circ$' and '$\ast$' markers show $|\phi|\cdot k$ for the circular duct at 400\,Hz and 800\,Hz. There is hardly any difference between the results for square and circular duct, in both ducts a frequency dependent damping the amplitude along the length of the duct can be clearly observed. Additionally, it can be seen that there is not only a damping effect, but also a frequency dependent offset in the curves. As there is no difference in the behavior between the results for the square and the circular duct, we will in the following concentrate on square ducts. 
\begin{table}[!h]
  \begin{center}
    \caption{Summary of the meshes used in Sec.~\ref{Sec:Surface}. The sidewalls $\Gamma$ are discretized using $N$ rectangular surface elements of size $h_1 \times h_2$, the closure at the ends of the duct $\Gamma_0$ and $\Gamma_L$ are discretized using square surface elements (square cross section) and triangular surface elements (circular cross section) with edge lengths $h_0$ and $h_L$. All edge lengths are given in m.}\label{Tab:SummarySurface}
  \begin{tabular}{ccccrccc}
    Nr. & Cross section & \multicolumn{3}{c}{$\Gamma$} & \multicolumn{3}{c}{$\Gamma_0, \Gamma_L$}\\
           & & $N$ & $h_1\um$ & $h_2\um$ & Elem. type & $N$ & $h_0\um,h_L\um$\\
    \hline
    1 & square & 1600 &0.04  &0.04 & quadrilateral & 25 & 0.04 \\
    2 & circular & 1600 & 0.04 & 0.035 & triangular & 64 & 0.035\\
  \end{tabular}
  \end{center}
\end{table}

\subsubsection{Effect of local grid refinement}\label{Sec:SurfLocal}
In \cite{Marburg18} the dependence of the damping effect on a local discretization of the duct was investigated. The duct was discretized with varying element sizes along the duct, and the field at the duct walls was investigated. In the following, we want to reproduce their results using the benchmark tube with $L = 3.4\tm.$

\begin{figure}[!h]
  \begin{center}
    \includegraphics[width=0.55\textwidth]{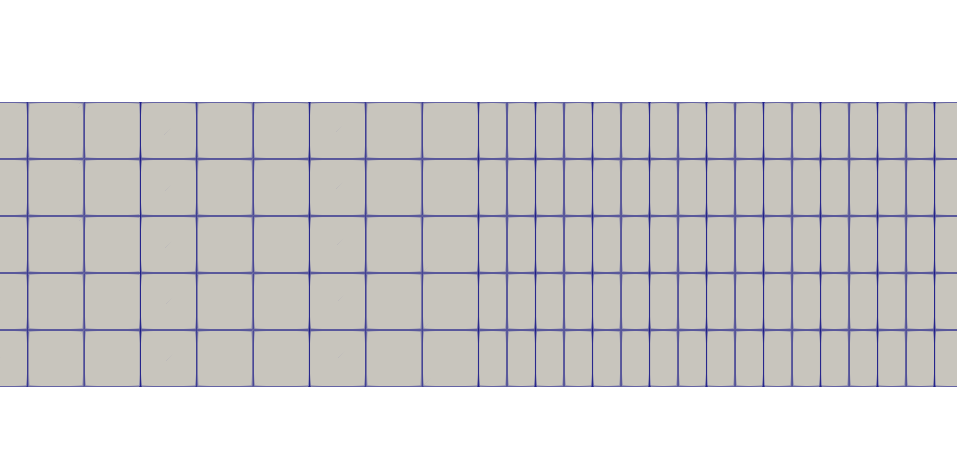}
    \caption{Mesh with element width $h_1 \approx 0.04\tm$ on the left part and $h_1 \approx 0.02\tm$ on the right part.}\label{Fig:FinevsCoarse}
  \end{center}
\end{figure}

The default height of each element used in the discretization of the sidewalls was $h_2 = 0.04\tm$. For the used frequency of $f = 800\,$Hz this corresponds to about 10 elements per wavelength. For the element width $h_1$  two different values $h_1 = 0.04\tm$ and $h_1 = 0.02\tm$ were used. First, the same $h_1$ was used uniformly along the duct, in a second experiment, the duct was divided into three parts where different $h_1$ were used. For the case $h_1 = (0.02, 0.04, 0.02)$ all rectangular elements in the first third of the duct with respect to its length  had a width of $h_1 \approx 0.02\tm$, in the second third of the duct they had a width of $h_1 \approx 0.04\tm$, and then again a width of $h_1 \approx 0.02\tm$ in the rest\footnote{As the duct was divided into 3 parts of equal length, and in each individual part each element had the same width, $h_1$ cannot be exactly $0.02\tm$ or $0.04\tm$.}. The definition for the case $h_1 = (0.04, 0.02, 0.04)$ is analogue. In Fig.~\ref{Fig:FinevsCoarse} the transition from a grid with element width $h_1 \approx 0.04\tm$ to a finer grid with element width $h_1 \approx 0.02\tm$ is depicted. The different meshes used in this subsection are summarized in Tab.~\ref{Tab:SurfLocal}.

As in Section~\ref{Sec:Surface} the particle velocity for elements at $x = 0$ was set to $v(\bx) = v_0 = 1\,$ms$^{-1}$, for all elements at $x = L = 3.4\tm$ an impedance of $Z = \rho c$ is assumed. For this problem the analytic solution has an absolute value of $|\phi_a(\bx)| = \frac{1}{k}\approx 0.067641$.  

\begin{figure}[!h]
  \begin{center}
    \includegraphics[width=0.85\textwidth]{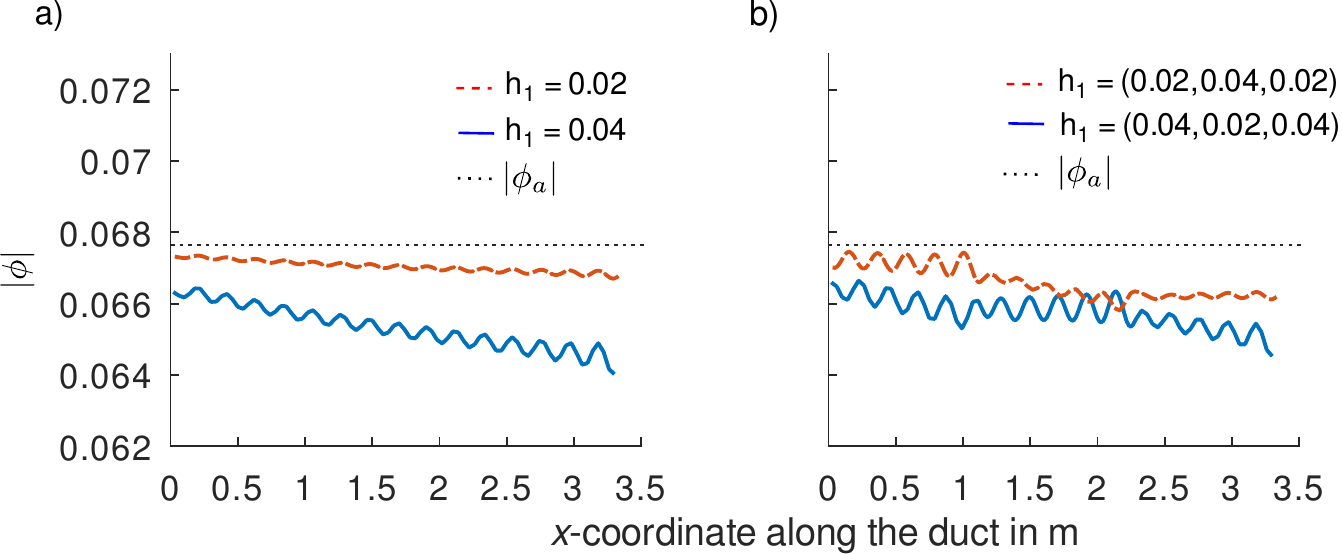}
  \end{center}
  \caption{Velocity potential along a sidewall of the square duct at 800\,Hz for different element sizes. a) $\phi$ is given for elements with width $h_1 = 0.02\tm$ (continuous line) and $h_1 = 0.04\tm$ (dashed line), respectively. b) The duct is divided into 3 parts with respect to its length, and different element sizes are used in different parts. The dotted lines show the constant value of the analytic solution $|\phi_a| = \frac{1}{k}$. All calculations have been done using surface elements.}\label{Fig:DiffDisc}
\end{figure}

In Fig.~\ref{Fig:DiffDisc} the velocity potential at the side of a square duct at a frequency $f = 800$\,Hz for different (local) element widths $h_1$  is depicted. In Fig.~\ref{Fig:DiffDisc}a, elements had a size of $h_1 \times h_2$ uniformly along the whole duct, $h_1 = 0.04\tm$ and $h_1 = 0.02\tm$, respectively. In Fig.~\ref{Fig:DiffDisc}b, different $h_1$ were used along the duct as described above.

As in \cite{Marburg18} we see that the damping of the amplitude is dependent on the discretization of the duct. But a smaller element width not only reduces the numerical damping, the ``level'' of the solution for a finer grid is also closer to the real solution. From  Fig.~\ref{Fig:DiffDisc}b) it becomes also clear that the effect of different $h_1$ is local. In parts of the duct with finer discretization, the damping effect is reduced, but as soon as the discretization becomes coarse, the damping effect gets bigger again.
\begin{table}[!h]
  \begin{center}
    \caption{Summary of the meshes used in Sec.~\ref{Sec:SurfLocal} and the figures they were used for. The sidewalls $\Gamma$ are discretized using rectangular surface elements of size $h_1 \times h_2$. The closures $\Gamma_0$ and $\Gamma_L$ are discretized using square surface elements with edge lengths $h_0$ and $h_L$. For meshes 3 and 4 the duct was split into three parts of equal length and the respective $h_1$ were used in each part. All lengths are given in m.}\label{Tab:SurfLocal}
    \begin{tabular}{ccccc}
      Nr. & $h_1\um$ & $h_2\um$ & $h_0\um,h_L\um$ & Figure\\
      \hline
      1 & 0.04 & 0.04 & 0.04 & Fig.~\ref{Fig:DiffDisc}a\\
      2 & 0.02 & 0.04 & 0.04 & Fig.~\ref{Fig:DiffDisc}a\\
      3 & 0.04, 0.02, 0.04 & 0.04 & 0.04 & Fig.~\ref{Fig:DiffDisc}b\\
      4 & 0.02, 0.04, 0.02 & 0.04 & 0.04 &  Fig.~\ref{Fig:DiffDisc}b
    \end{tabular}
  \end{center}
\end{table}
\subsubsection{Field inside the duct}\label{Sec:InsideSurf}
\begin{figure}[!h]
  \begin{center}
    \includegraphics[width=0.90\textwidth]{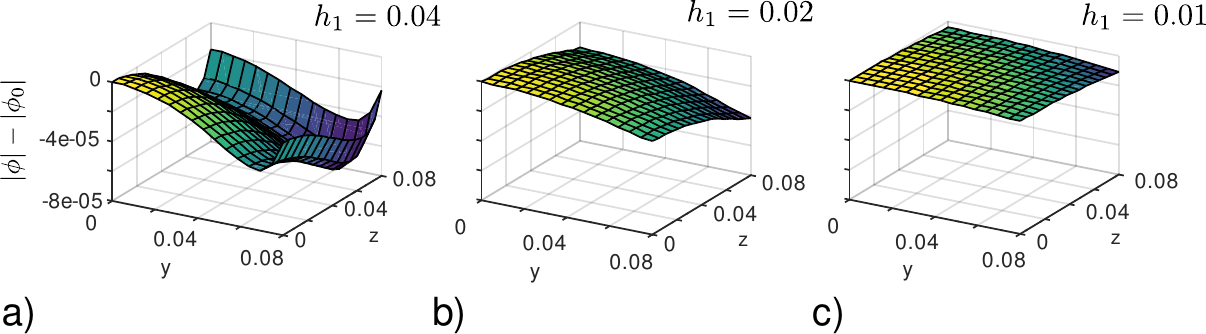}\\[12pt]
    \includegraphics[width=0.90\textwidth]{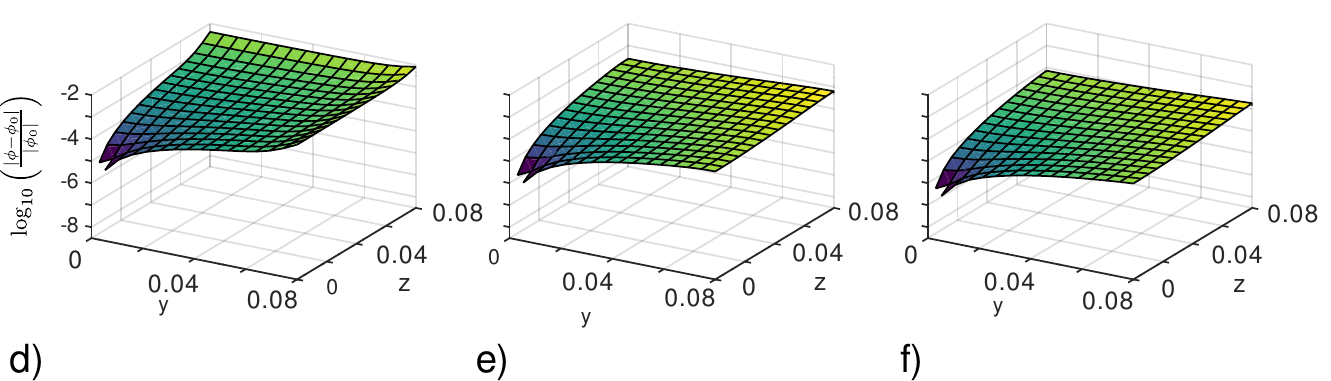}
  \end{center}
  \caption{Difference of the velocity potential $\phi$ in the plane normal to the duct axis inside the duct at $x_0 = \frac{L}{3}$ and $\phi_0 = \phi(x_0,0,0)$. In the upper row the difference of the absolute values is depicted, in the lower row the logarithm of the relative distance between the $\phi$ and $\phi_0$. In the first column, an element width of $h_1 = 0.04\tm$ is used, in the second column $h_1 = 0.02\tm$, and in the third columns $h_1 = 0.01\tm$. As the field is symmetric with respect to the $y$ and $z$ axis, only a quarter of the plane is displayed. All calculations have been done using surface elements.}\label{Fig:InsideSurf}
\end{figure}

As for this example the field inside the duct is given by a single plane wave, one would expect the absolute value of the velocity potential to be constant in each plane normal to the duct axis. However, in the simulation $|\phi|$ gets smaller towards the duct sidewalls, as can be seen in Fig.~\ref{Fig:InsideSurf}, where the field at a plane inside the duct normal to the duct axis for a frequency of $f =  800\,$Hz is shown. 

In the first row the shifted absolute values of the velocity potential $\phi$ at the evaluation points at $x_0 = \frac{L}{3}$ are depicted. For better visualization and because of the damping effect, the absolute value $|\phi_0| = |\phi(x_0,0,0)|$ at the middle line of the duct was subtracted from the velocity potential at the evaluation points $|\phi(x_0,y_i,z_i)|$, where $y_i$ and $z_i$ ranged from $0.0\tm$ to $0.08\tm$ as described at the beginning of Section~\ref{Sec:ExperimentI}. The difference of the absolute values was chosen to better illustrate the deviation of the BEM solution from the plane wave solution. 
From left to right the width of each BE-element was set to $h_1 = 0.04\tm$, $h_1 = 0.02\tm$, and $h_1 = 0.01\tm$, respectively. The height of each element was always kept at $h_2 = 0.04\tm$. The meshes used in this section have been summarized in Tab~\ref{Tab:SurfInside}.

In the second row of Fig.~\ref{Fig:InsideSurf} the relative distance $\frac{|\phi - \phi_0|}{|\phi_0|}$ of the velocity potentials is depicted. It can be observed that this relative distance becomes smaller  with smaller $h_1$, thus the BEM calculations seem to converge towards the ``correct'' plane wave solution.
\begin{table}[!h]
  \begin{center}
    \caption{Summary of the meshes used in Sec.~\ref{Sec:InsideSurf} and the figures they were used for. The sidewalls $\Gamma$ are discretized using rectangular surface elements of size $h_1 \times h_2$. The closures $\Gamma_0$ and $\Gamma_L$ are discretized using square surface elements with edge lengths $h_0$ and $h_L$. All lengths are given in m.}\label{Tab:SurfInside}
    \begin{tabular}{ccccc}
      Nr. & $h_1\um$ & $h_2\um$ & $h_0\um,h_L\um$ & Figure\\
      \hline
      1 & 0.04 & 0.04 & 0.04 & Figs.~\ref{Fig:InsideSurf}a,\ref{Fig:InsideSurf}d\\
      2 & 0.02 & 0.04 & 0.04 & Figs.~\ref{Fig:InsideSurf}b,\ref{Fig:InsideSurf}e\\
      3 & 0.01 & 0.04 & 0.04 & Figs.~\ref{Fig:InsideSurf}c,\ref{Fig:InsideSurf}f
    \end{tabular}
  \end{center}
\end{table}
\subsection{Duct with thin elements}\label{Sec:Thin}
\subsubsection{Field along the duct}\label{Sec:ThinAlong}
\begin{figure}[!h]
  \begin{center}
    \includegraphics[width=0.95\textwidth]{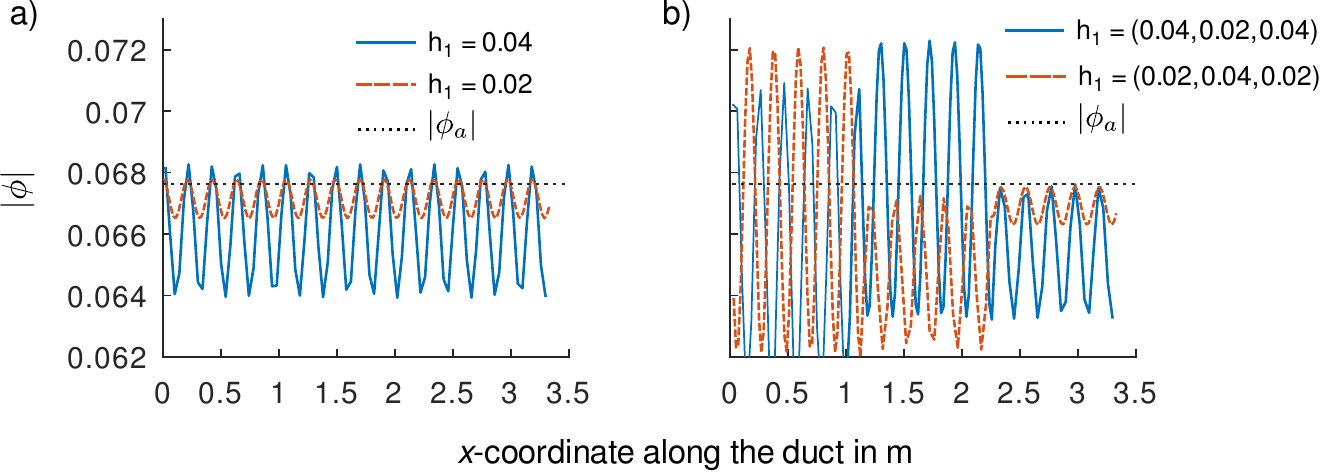}
  \end{center}
  \caption{Velocity potential along one sidewall of the square duct at 800\,Hz for different element widths using thin elements. a) $\phi$ is given for elements with width $h_1 = 0.02\tm$ (continuous line) and $h_1 = 0.04\tm$ (dashed line), respectively. b) The duct in divided into 3 parts with respect to its length, and different element widths are used in different parts. The dotted line shows the constant value of the analytic solution $|\phi_a|$. All calculations have been done using thin elements.}\label{Fig:Midfacesol}
\end{figure}

Using thin elements instead of regular surface elements results in a slightly different behavior of the velocity potential. First, there is a greater variation of the numerical solution in comparison with the solution for surface elements, and second, a damping of the amplitude is not apparent, see  Fig.~\ref{Fig:Midfacesol}, where the same setup as in Section~\ref{Sec:SurfLocal} was used. The meshes used in this section have been summarized at the end of the section in Tab.~\ref{Tab:Thin}.
The variation can be reduced by a great amount by reducing the width of the elements along the \emph{whole} duct, which is clearly visible in Fig.~\ref{Fig:Midfacesol}a. The continuous line depicts the absolute value of the velocity potential along a line on the side of the duct. The sidewalls of the duct are discretized with elements with  $h_1 = 0.04\tm$ along the duct. For the dashed line the width of each element on the sidewalls was $h_1 = 0.02$\,m. In both cases the height of the elements was fixed at $h_2 = 0.04$\,m. 

As for surface elements, there is a local behavior with respect to different element widths, but the direct relationship between  element width and error is not that obvious anymore. Using different $h_1$ in different parts along the duct still results in local changes, but this one-to-one relation between $h_1$ and the size of the error as observed for surface elements is not given anymore. A small $h_1$ in one part of the duct does not necessarily mean, that the variation of the solution in this part is small too. 

\begin{figure}[!h]
  \begin{center}
    \includegraphics[width=0.85\textwidth]{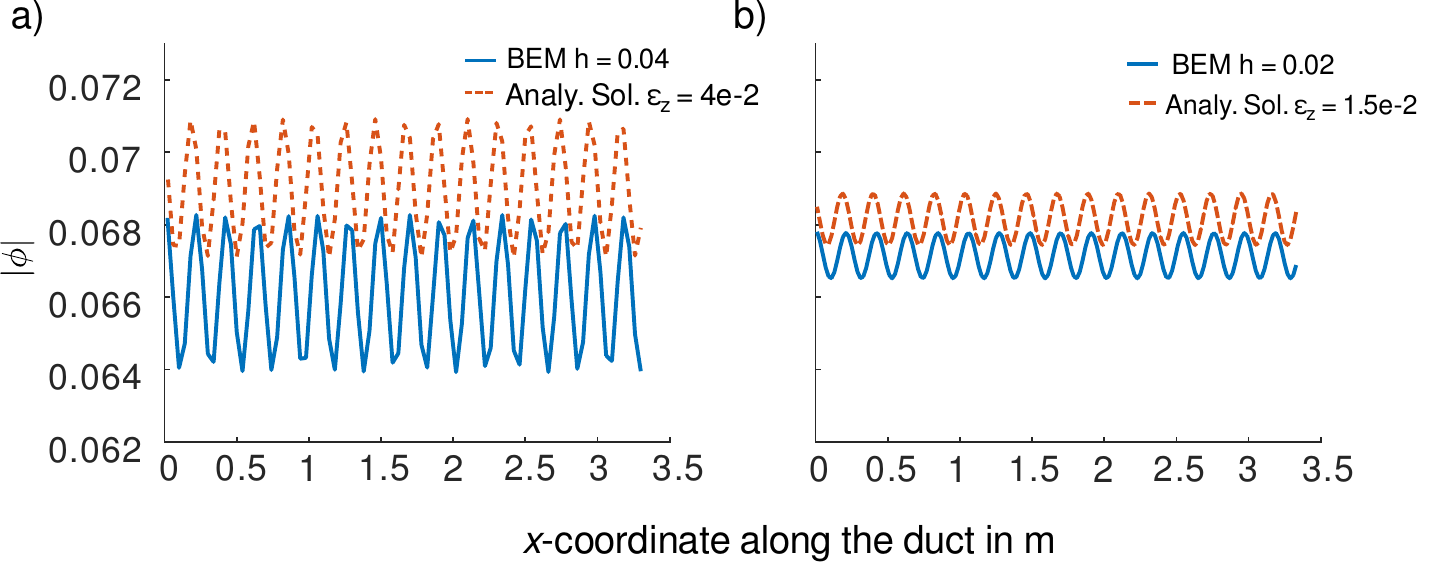} 
    \caption{BEM solution using thin elements (continuous line) and the analytic solution with a perturbed impedance $Z_p = (1 + \varepsilon_Z)\rho c$ at $x = L$ (dashed line) for element widths of a) $h_1 = 0.04\tm$ and b) $h_1 = 0.02\tm$. All calculations have been done using thin elements.}\label{Fig:BEMvsAnaly}
  \end{center}
\end{figure}

A close look at Fig.~\ref{Fig:Midfacesol}b suggests that the variation of the amplitude in the third part of the duct correlates with the width of the  elements close to $x = L$. One explanation for that behavior may be that the impedance $Z(L) = \rho c$  used at the closure of the duct at $x = L$ gets perturbed, resulting in an additional plane wave traveling in the opposite direction. In Fig.~\ref{Fig:BEMvsAnaly}  BEM solutions using thin elements with $h_1 = 0.04\tm$ (left side) and $h_1 = 0.02\tm$ (right side) are compared with the analytic 1D-solution for ducts using Eq.~(\ref{Equ:Helmholtz1D}), where at $x = L$ a perturbed impedance $Z_p = ( 1 + \varepsilon_Z) \rho c$ with $\varepsilon_Z = 4\cdot 10^{-2}$ and $\varepsilon_Z = 1.5 \cdot 10^{-2}$ was assumed. These solutions are depicted using the dashed lines. As one can see, the numerical and the analytic solutions have about the same level of variation.
\begin{table}[!h]
  \begin{center}
      \caption{Summary of the meshes used in Sec.~\ref{Sec:ThinAlong} and the figures they were used for. The sidewalls $\Gamma$ are discretized using rectangular thin elements of size $h_1 \times h_2$. The closures $\Gamma_0$ and $\Gamma_L$ are discretized using square thin elements with edge lengths $h_0$ and $h_L$. For meshes 3 and 4 the duct was split into three parts of equal length and the respective $h_1$ were used in each part. All lengths are given in m.}\label{Tab:Thin}
    \begin{tabular}{ccccc}
      Nr. & $h_1\um$ & $h_2\um$ & $h_0\um,h_L\um$ & Figure\\
      \hline
      1 & 0.04 & 0.04 & 0.04 & Fig.~\ref{Fig:Midfacesol}a, \ref{Fig:BEMvsAnaly}a\\
      2 & 0.02 & 0.04 & 0.04 & Fig.~\ref{Fig:Midfacesol}a, \ref{Fig:BEMvsAnaly}b\\
      3 & 0.04, 0.02, 0.04 & 0.04 & 0.04 & Fig.~\ref{Fig:Midfacesol}b\\
      4 & 0.02, 0.04, 0.02 & 0.04 & 0.04 & Fig.~\ref{Fig:Midfacesol}b
    \end{tabular}
  \end{center}
\end{table}

\subsubsection{Influence of the discretization of the duct closures}\label{Sec:ThinClosure}
There is a connection between the accuracy of the numerical solutions and the discretizations of the duct closures at $x = 0$ and $x = L$, that can only be found when using thin elements. In short, results for thin boundary elements get closer to the analytic solution if the closures $\Gamma_0$ and $\Gamma_L$ at $x = 0$ and $x = L$ are discretized using \emph{less} elements.

\begin{figure}[!h]
  \begin{center}
    \includegraphics[width=0.9\textwidth]{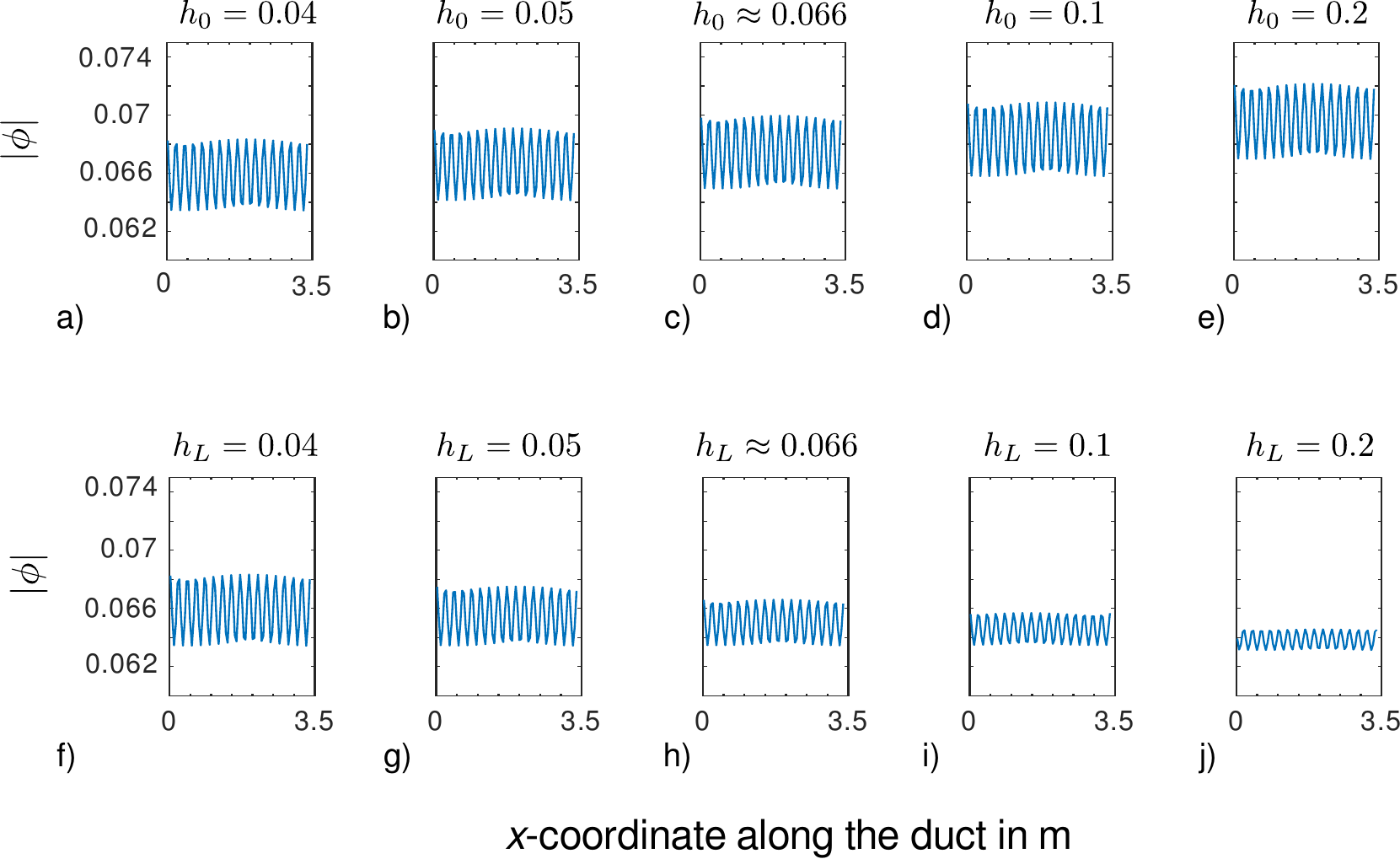}
  \end{center}
  \caption{Comparison of the effect of different discretizations of the closures $\Gamma_0$ and $\Gamma_L$ at $x = 0$ and $x = L$, respectively. In the first row the discretization of the closure $\Gamma_0$ gets coarser from left to right, in the second row the discretization of the closure $\Gamma_L$. $h_0$ and $h_L$ denote the edge lengths of the square elements used in $\Gamma_0$ and $\Gamma_L$, respectively. In all calculations thin elements have been used.}\label{Fig:VeloPotDiffClosure}
\end{figure}

In Fig.~\ref{Fig:VeloPotDiffClosure} the effect of different discretizations of the duct closures is shown. In the first row the effect of different discretization of the closure $\Gamma_0$ at $x = 0$, where the velocity condition is applied, is depicted. The discretization gets \emph{coarser} from left to right, starting with $5\times 5$ elements with an element size of $h_0 \times h_0 = 0.04\tm \times 0.04\tm$ in Fig.~\ref{Fig:VeloPotDiffClosure}a to just 1 element with an element size of $h_0 \times h_0 = 0.2\tm \times 0.2\tm$ in Fig.~\ref{Fig:VeloPotDiffClosure}e. The discretization of the right end of the duct at $x = L$ is always kept to $5\times 5$ elements with a size of $h_L\times h_L = 0.04\tm \times 0.04\tm$. 

In the second row (Figs.~\ref{Fig:VeloPotDiffClosure}f to \ref{Fig:VeloPotDiffClosure}g), $|\phi|$ is shown for different element sizes $h_L \times h_L$ at $x = L$ and fixed element size $h_0 \times h_0 = 0.04\tm \times 0.04\tm$ at $x = 0$. A coarse discretization of the closure at $x = L$ reduces the variation in the velocity potential noticeable. The meshes used are summarized in Tab.~\ref{Tab:ThinClosure}.

This behavior may seem odd at a first glance, but some remarks need to be added. First, the fact that only one element is used in the discretization of the closure, does \emph{not} mean, that the quadrature over the closures is not accurate. In our code an adaptive quadrature scheme similar to the method in \cite{LacWat76} is used, that is based on an apriori error estimator to determine the order of the Gauss quadrature and that additionally subdivides elements if the estimated Gauss order is bigger than a fixed maximum order. Using just one element for the closures in combination with the collocation approach with constant elements has the effect that the numerical solution is forced to be constant over the closures, which in turn suggests, that the solution is closer to the constant field of a plane wave and less numerical perturbation of the impedance boundary condition is introduced. A similar argument can be used for the discretization of the closure at $x = 0$ that has an influence on the overall level of the velocity potential.

\begin{figure}[!h]
  \begin{center}
     \includegraphics[width=0.75\textwidth]{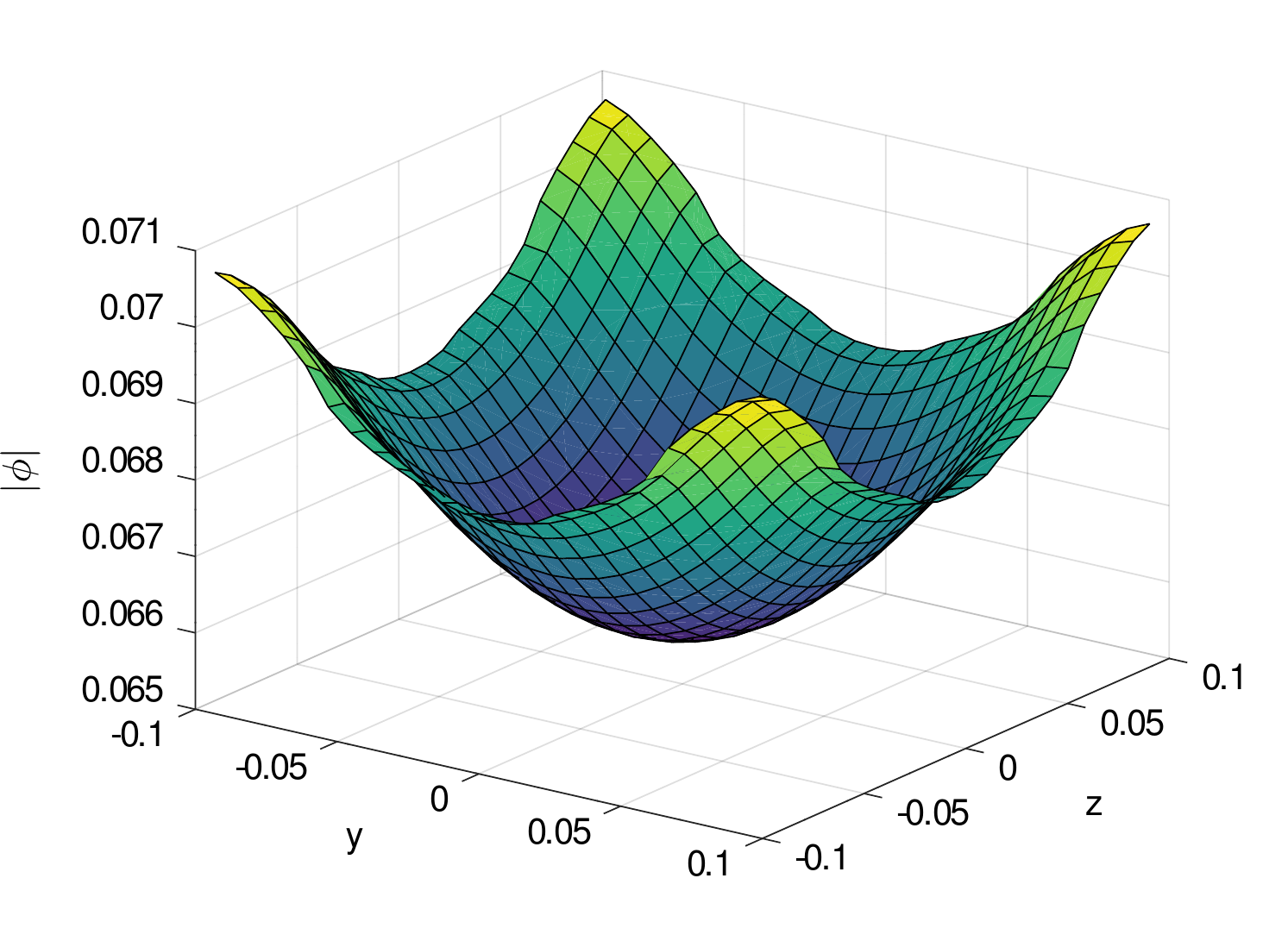}
  \end{center}
  \caption{Velocity potential at the closure at $x = L$. The sidewalls of the duct were discretized using thin elements with default element size $h_1 \times h_2 = 0.04\tm \times 0.04\tm$, the closures were discretized using $25\times 25$ thin elements of size $h_0 \times h_0 = h_L \times h_L = 0.008\tm \times 0.008\tm$.}\label{Fig:PotEnd}
\end{figure}

In Fig.~\ref{Fig:PotEnd} the velocity potential  at $f = 800$\,Hz  on the closure $\Gamma_L$ for a very fine discretization of both closures of the duct is displayed. Each sidewall is discretized using the standard element width and height $h_1 \times h_2 = 0.04\tm \times 0.04\tm$, the closures of the duct, on the other hand, are discretized using $25\times 25$ elements of size $0.008\tm \times 0.008\tm$. The absolute value of $\phi$ at $x = L$ for this discretization ranges from 0.065682 to 0.070702, this range becomes smaller for coarser discretizations of the closures, if just one element is used the absolute value of the velocity potential is about $0.067741$ and the difference to the absolute value the analytic solution is smaller than $10^{-4}$.
For surface elements the field shows similar behavior, however, minimum and maximum values of the absolute values lie between approximately 0.064 and 0.066, showing less variation but the numerical damping effect. The fact that the calculated field with a fine discretization of the closure is not constant over the cross section may explain that with a coarser discretization of $\Gamma_L$ the variation in the velocity potential becomes smaller. Quoting \cite{Fahnline08}, who investigated BEM with surface elements applied to ducts: ``However, in broad terms, the solution accuracy is still determined by how well the specified boundary condition is satisfied''. This seems to be the case for surface as well as thin elements.

\begin{table}[!h]
  \caption{Summary of the meshes used in Sec.~\ref{Sec:ThinClosure} and the figures they were used for. The sidewalls $\Gamma$ are discretized using rectangular thin elements of size $h_1 \times h_2$. The closures $\Gamma_0$ and $\Gamma_L$ are discretized using square thin elements with edge lengths $h_0$ and $h_L$. All lengths are given in m.}\label{Tab:ThinClosure}
  \begin{center}
    \begin{tabular}{rccrrl}
      Nr. & $h_1\um$ & $h_2\um$ & $h_0\um$ & $h_L\um$ & Figure\\
      \hline
      1 & 0.04 & 0.04 & 0.04  & 0.04 & Fig.~\ref{Fig:VeloPotDiffClosure}a\\
      2 & 0.04 & 0.04 & 0.05  & 0.04 & Fig.~\ref{Fig:VeloPotDiffClosure}b\\
      3 & 0.04 & 0.04 & 0.066 & 0.04 & Fig.~\ref{Fig:VeloPotDiffClosure}c\\
      4 & 0.04 & 0.04 & 0.1   & 0.04 & Fig.~\ref{Fig:VeloPotDiffClosure}d\\
      5 & 0.04 & 0.04 & 0.2   & 0.04 & Fig.~\ref{Fig:VeloPotDiffClosure}e\\
      7 & 0.04 & 0.04 & 0.04 & 0.04  & Fig.~\ref{Fig:VeloPotDiffClosure}f\\
      8 & 0.04 & 0.04 & 0.04 & 0.05  & Fig.~\ref{Fig:VeloPotDiffClosure}g\\
      9 & 0.04 & 0.04 & 0.04 & 0.066 & Fig.~\ref{Fig:VeloPotDiffClosure}h\\
      10 & 0.04 & 0.04& 0.04  & 0.1  & Fig.~\ref{Fig:VeloPotDiffClosure}i\\
      11 & 0.04 & 0.04& 0.04  & 0.2  & Fig.~\ref{Fig:VeloPotDiffClosure}j\\
      12 & 0.04 & 0.04 & 0.008 & 0.008 & Fig.~\ref{Fig:PotEnd}
    \end{tabular}
  \end{center}
\end{table}

\subsubsection{Field inside the duct}\label{Sec:ThinInside}
\begin{figure}[!h]
  \begin{center}
  \includegraphics[width=0.9\textwidth]{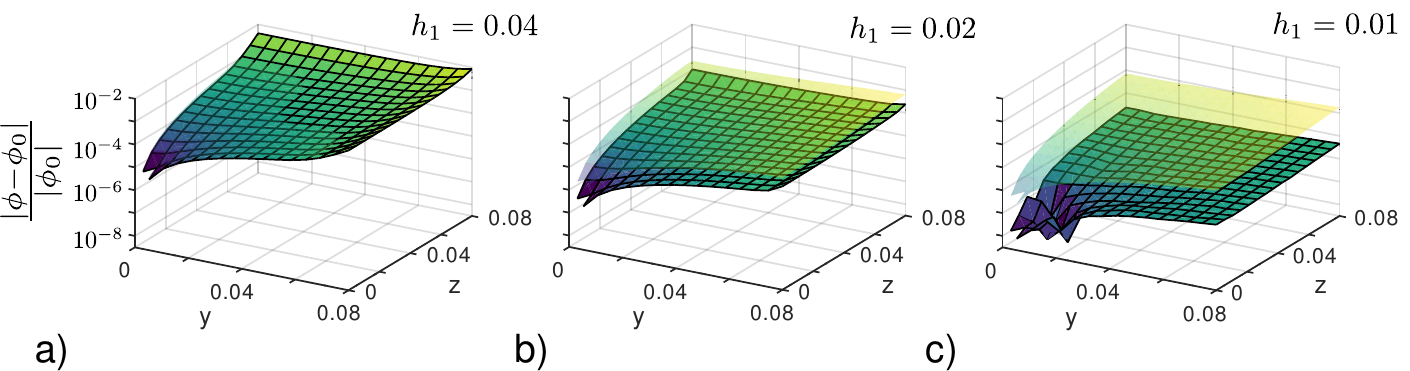}
  \caption{Relative distance between the velocity potential $\phi$ inside the duct on a plane normal to the duct axis at $x_0 = \frac{L}{3}$ and $\phi_0 = \phi(x_0,0,0)$ as a function of element width $h_1$. a) $h_1 = 0.04\tm$, b) $h_1 = 0.02\tm$, c) $h_1 = 0.01\tm$. For comparison with the results for surface elements, the transparent patches show the distance for surface elements, the other patches the calculations using thin elements.}\label{Fig:InsideMid}
  \end{center}
\end{figure}

Using the same setup as in Sec.~\ref{Sec:InsideSurf}, the field inside the duct for thin elements behaves slightly differently to the case when surface elements are used. First, when looking at middle line of the duct, the velocity potential behaves like the field on the sides of the duct. There is no damping effect visible, but a greater variation of the amplitude of the velocity potential. Second, the amplitude of the acoustic field  normal to the duct axis, is, in general, closer to a constant plane wave solution than the field calculated using surface elements, as can be observed in Fig.~\ref{Fig:InsideMid}. The convergence towards a constant value is much quicker as a function of element width $h_1$ for the thin elements. For an easier comparison with the surface element solution the transparent patches in Fig.~\ref{Fig:InsideMid} depict the solutions for BEM with surface elements and the same $h_1$. The meshes used in this section are summarized in Tab.~\ref{Tab:ThinInside}.
\begin{table}[!h]
  \begin{center}
    \caption{Summary of the meshes used in Sec.~\ref{Sec:ThinInside} and the figures they were used for. The sidewalls $\Gamma$ are discretized using rectangular thin elements of size $h_1 \times h_2$. The closures $\Gamma_0$ and $\Gamma_L$ are discretized using square thin elements with edge lengths $h_0$ and $h_L$. All lengths are given in m.}\label{Tab:ThinInside}
    \begin{tabular}{ccccc}
      Nr. & $h_1\um$ & $h_2\um$ & $h_0\um,h_L\um$ & Figure\\
      \hline
      1 & 0.04 & 0.04 & 0.04 & Figs.~\ref{Fig:InsideMid}a\\
      2 & 0.02 & 0.04 & 0.04 & Figs.~\ref{Fig:InsideMid}b\\
      3 & 0.01 & 0.04 & 0.04 & Figs.~\ref{Fig:InsideMid}c
    \end{tabular}
  \end{center}
\end{table}
\subsection{Error along the duct}\label{Sec:ErrorAlong}
To better illustrate the behavior of the solutions for surface and thin elements the length of the duct is changed to 10\,m just in this section. Each boundary element used in the discretization of the duct has a height of $h_2 = 0.04$\,m, for its width $h_1 = 0.04$\,m, $h_1 = 0.02$\,m, and $h_1 = 0.01$\,m  are used, respectively. The values of the BEM solution $\phi$ are evaluated at collocation points along a line on the surface of a sidewall at $\bx = (x,y,z)$ with $x$\ lies between $0\tm$ and $10\tm$,  $y = -0.1\tm$, and $z = 0\tm$. The plane wave solution $\phi_a = \frac{1}{\I k} \e^{-\I k x}$ is evaluated at the same line, the relative error is defined by $\varepsilon_r(x) = \frac{|\phi(x) - \phi_a(x)|}{\phi_a(x)}$.

\begin{figure}[!h]
  \begin{center}
    \includegraphics[width=0.40\textwidth]{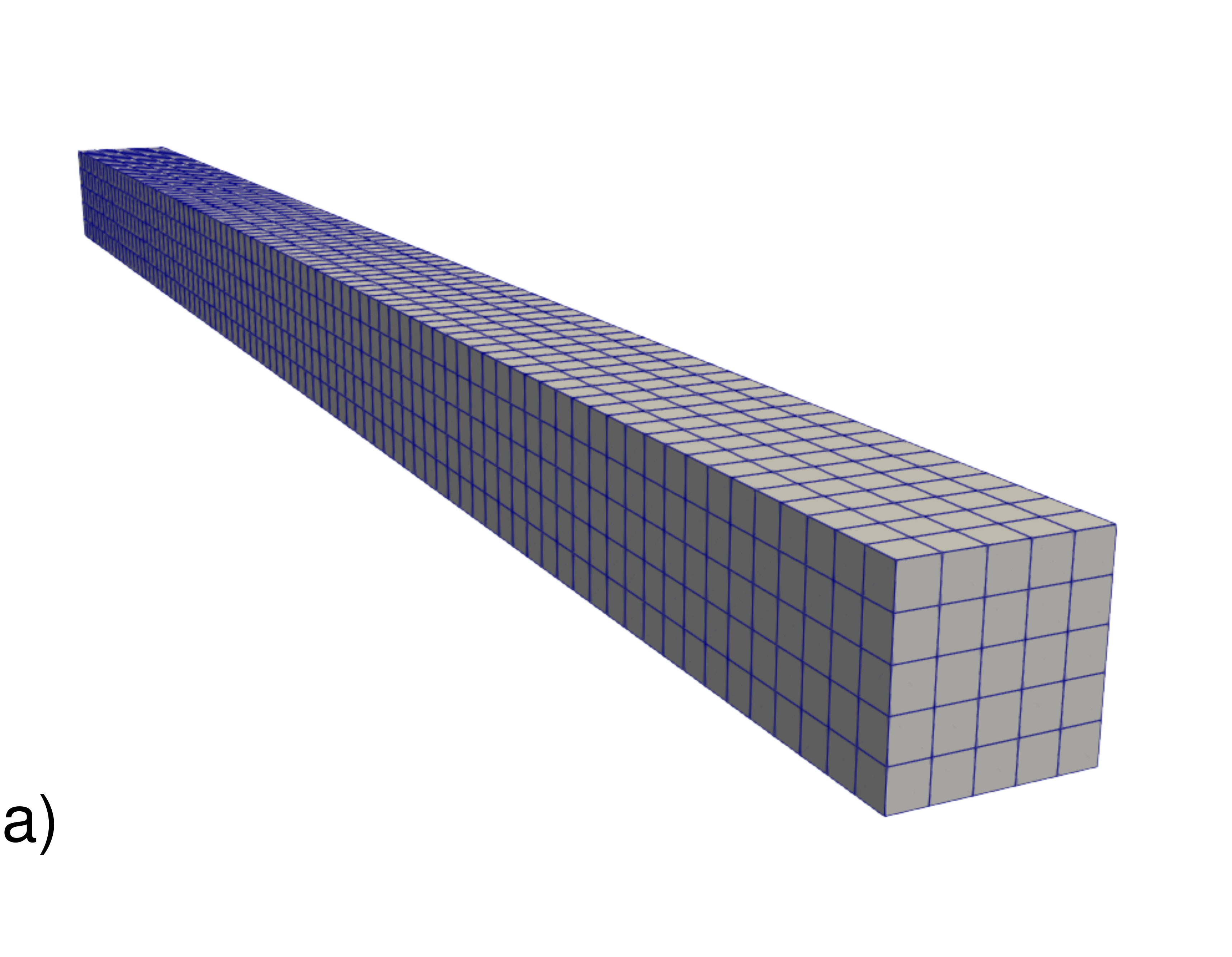}
    \includegraphics[width=0.40\textwidth]{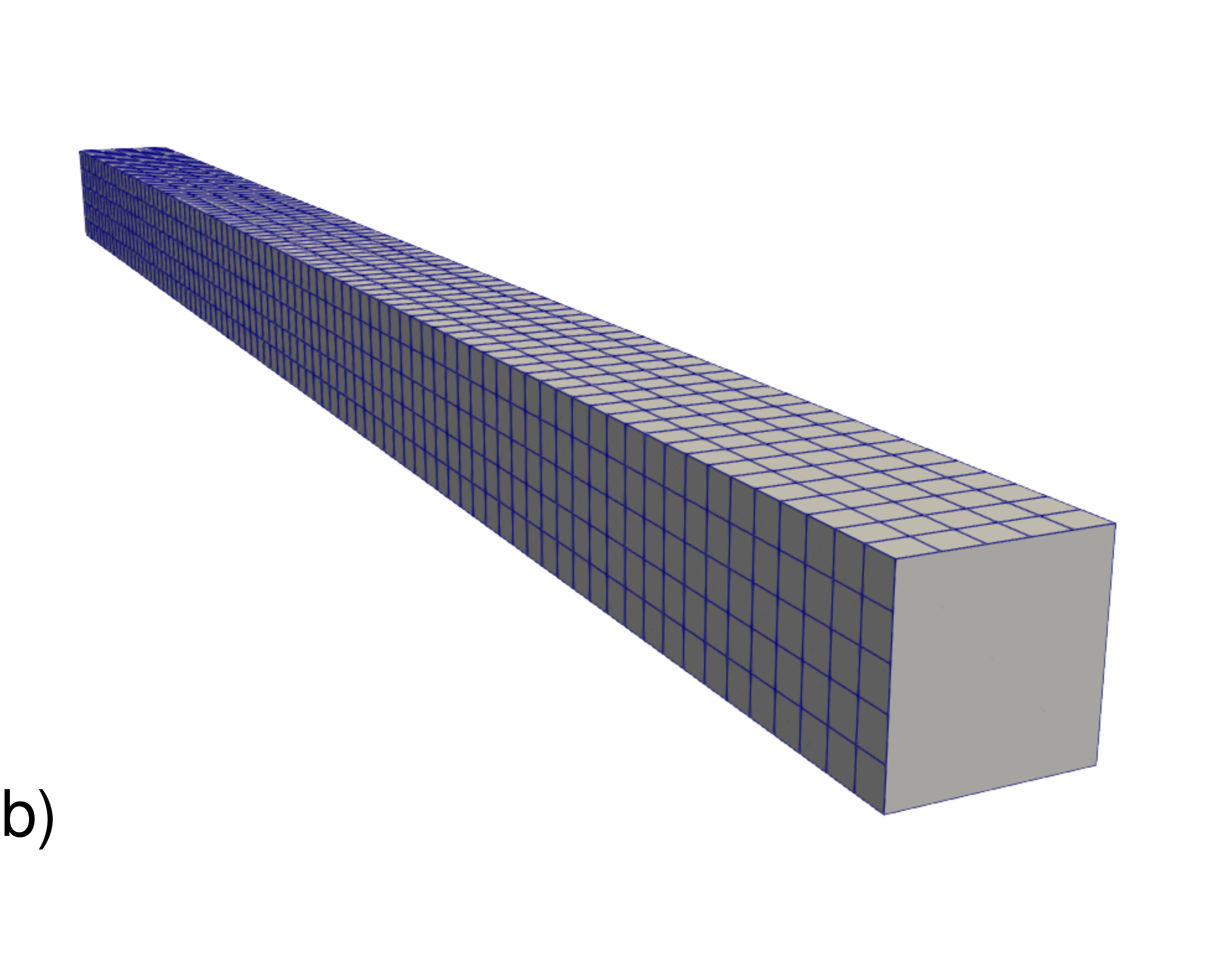}
  \end{center}
  \caption{Comparison between the different meshes used in this section. a) Surf$_5$ and Thin$_5$ mesh, b) Thin$_1$ mesh.}\label{Fig:TheThins}
\end{figure}

As it became clear in Sec.~\ref{Sec:Thin}, the discretization of the duct closures has an effect on overall level and variation of the solution for thin elements. Thus, we will in the following concentrate on three types of discretizations (the different meshes are summarized in Tab.~\ref{Tab:ErrorAlong} and depicted in Fig.~\ref{Fig:TheThins}):
\begin{description}
\item[Thin$_5$:] The sides of the duct are discretized using thin boundary elements with size $h_1 \times h_2$, the closures $\Gamma_0$ and $\Gamma_L$ at $x = 0$ and $x = L$ are discretized by using $5\times 5$ thin quadrilateral elements, thus each element has the size $h_0 \times h_0 = h_L \times h_L = 0.04\,\text{m} \times 0.04\,\text{m}$.
\item[Thin$_1$:] The sides of the duct are discretized using thin boundary elements with size $h_1 \times h_2$, the closures $\Gamma_0$ and $\Gamma_L$ at $x = 0$ and $x = L$ are discretized with just 1 element, thus, $h_0 \times h_0 = h_L \times h_L = 0.2\,\text{m} \times 0.2\,\text{m}$.
\item[Surf$_5$:]The sides of the duct are discretized using regular surface boundary elements with size $h_1 \times h_2$, the closures $\Gamma_0$ and $\Gamma_L$ at $x = 0$ and $x = L$ are discretized by using $5\times 5$ thin quadrilateral elements, thus each element has the size $h_0 \times h_0 = h_L \times h_L = 0.04\tm \times 0.04\tm$. 
\end{description}

\begin{figure}[!h]
  \includegraphics[width=0.95\textwidth]{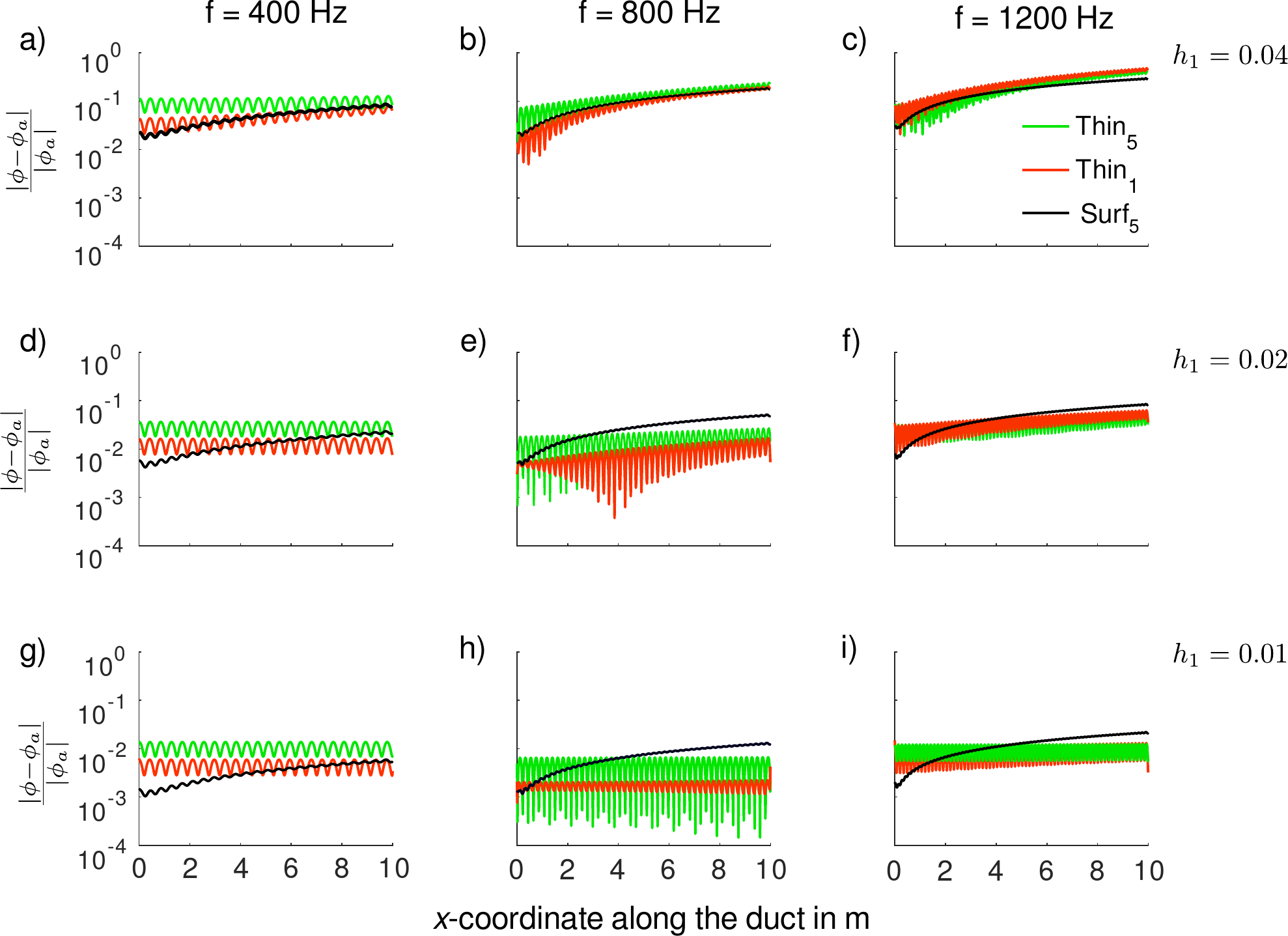}
  \caption{Relative distance $\frac{|\phi_a - \phi|}{|\phi_a|}$ of the BEM calculations $\phi$ and the analytic solution $\phi_a$ along a line on the duct wall for the different numerical calculations using thin and surface elements. For the graphs in the first row (a,b and c) the width of each boundary element was $h_1 = 0.04\tm$, in the second row (e,f, and g) $h_1 = 0.02\tm$, and in the third row (g,h, and i) $h_1 = 0.01\tm$. In all cases the height of each boundary element was fixed to $h_2 = 0.04\tm$. In the first column the errors are calculated for a frequency $f = 400$\,Hz, in the second column for $f = 800$\,Hz, and in the third column for $f = 1200$\,Hz. The green (light gray) curves show the error curves for a discretization using Thin$_5$ elements, the red (dark gray) lines for a Thin$_1$ discretization, and the black lines the errors for calculations using surface elements.}\label{Fig:RelErrorDiffn1}
\end{figure}

\begin{figure}[!h]
  \begin{center}
    \begin{center}
      \includegraphics[width=0.48\textwidth]{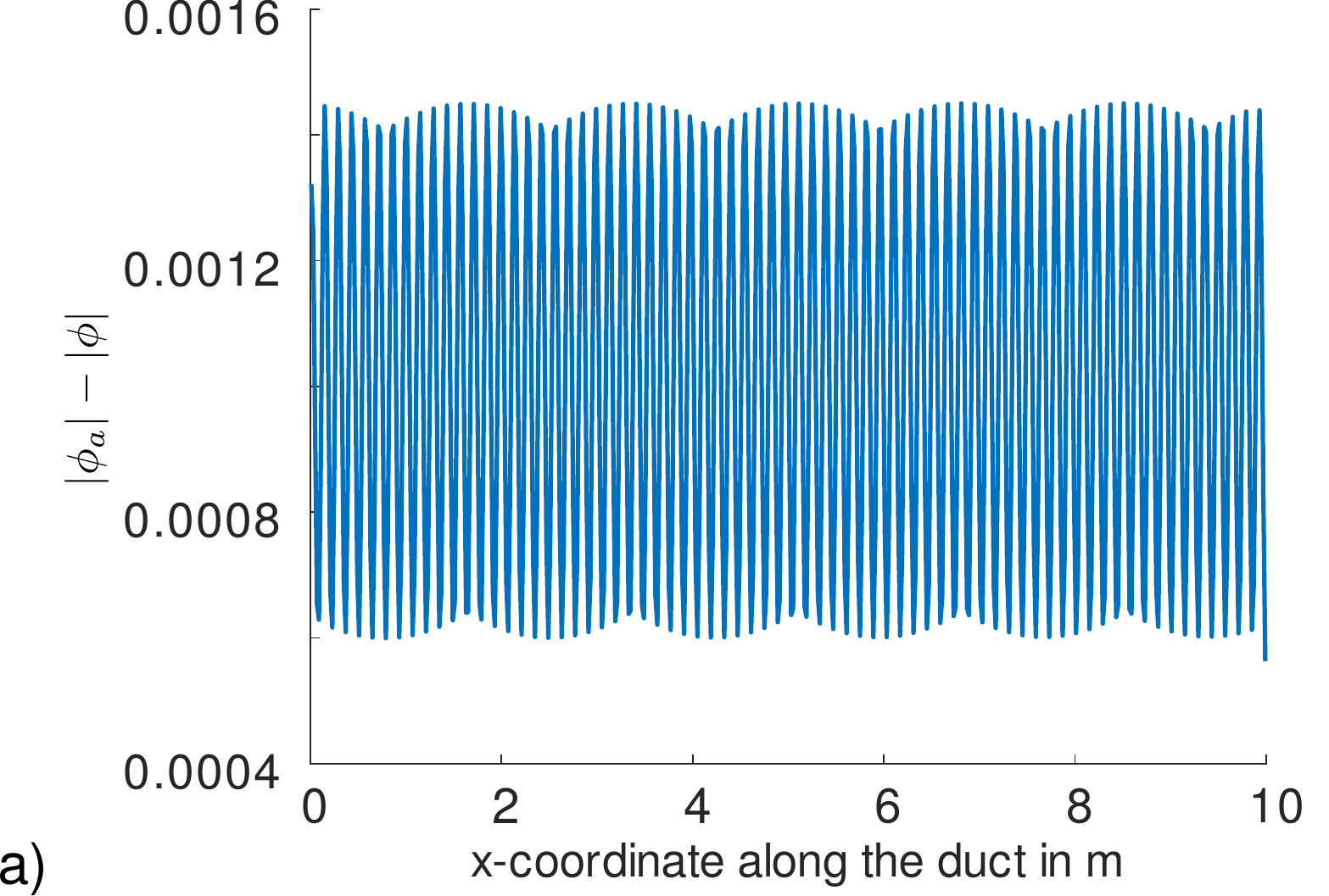}
      \includegraphics[width=0.48\textwidth]{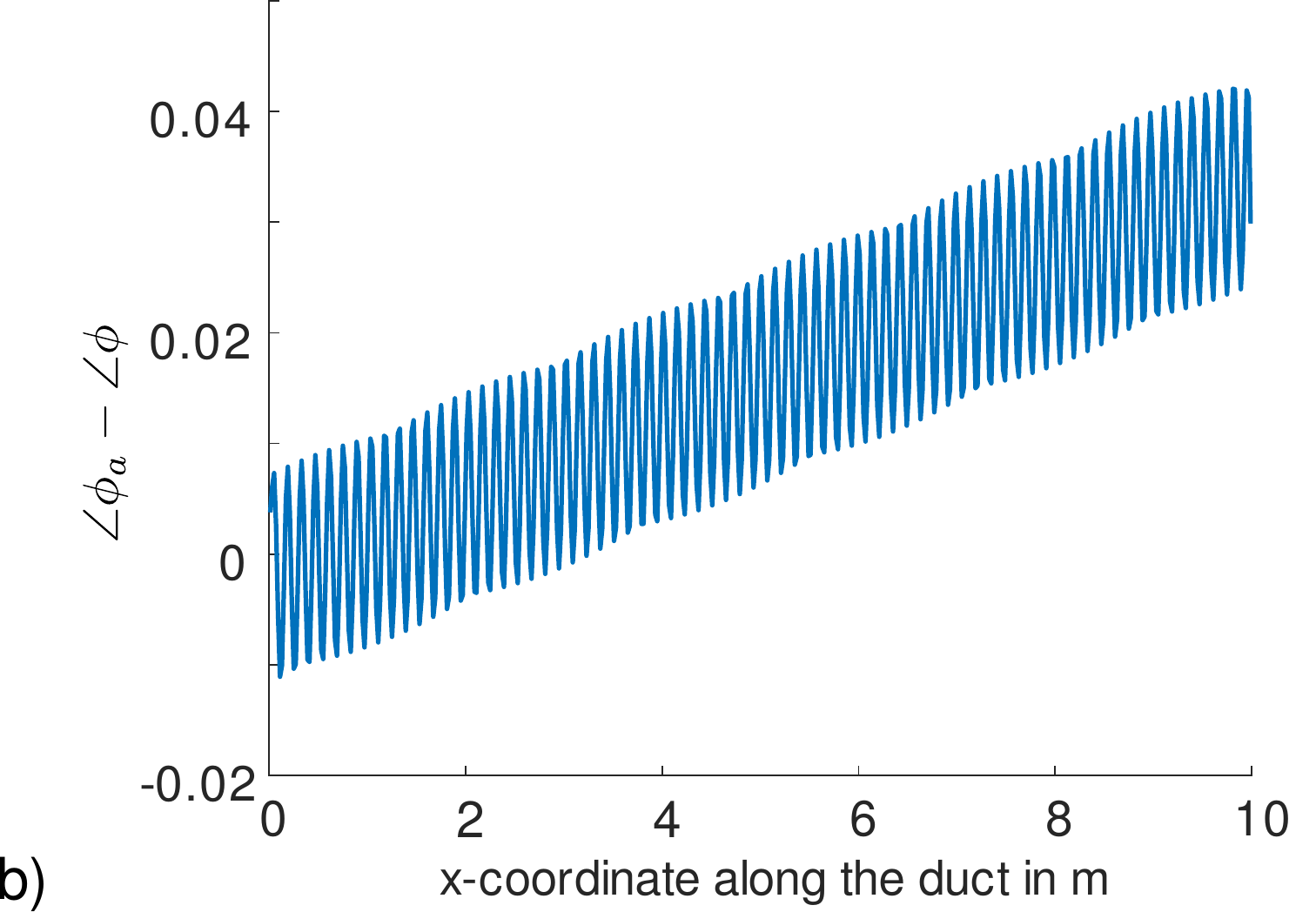}
    \end{center}
    \caption{a) Difference of the absolute values  between the analytic solution $\phi_a$ and the BEM solution $\phi$. b) Difference in rad of the polar angles between the analytic solution $\phi_a$ and the BEM solution $\phi$.}\label{Fig:Anglealong}
  \end{center}
\end{figure}

We want to point out, that there \emph{is} an accumulation of the  numerical error along the length of the duct when using thin elements, see, e.g., Fig~\ref{Fig:RelErrorDiffn1}. In this figure the relative error $\frac{|\phi - \phi_a|}{|\phi_a|}$ of the numerically calculated velocity potential along the duct using BEM with thin and surface elements is depicted, where $\phi$ and $\phi_a$ are the  the numerical and the analytic solutions, respectively. For the curves in the first row, $h_1 = 0.04\tm$, in the second row $h_1$ was set to $h_1 = 0.02\tm$, and in the third row $h_1 = 0.01\tm$, $h_2$ was always kept at $h_2 = 0.04\tm$. Each column in Fig.~\ref{Fig:RelErrorDiffn1} shows the error for different frequencies $f$. In the first column $f = 400$\,Hz, in the second $f = 800$\,Hz, and in the third $f = 1200$\,Hz. The green (light gray) curves show the errors for Thin$_5$ elements, the red (dark gray) lines for Thin$_1$ elements, and the black lines for surface elements.
In the previous sections it became clear that for thin elements the amplitude of $\phi$ seems to have a constant mean value along the duct. The visible rise in the error in Fig.~\ref{Fig:RelErrorDiffn1} also for thin elements suggests that this rise in error is mainly caused by shifts in the phase, which can also be observed when looking at Fig.~\ref{Fig:Anglealong}. In Fig.~\ref{Fig:Anglealong}a the difference of the absolute values of the analytic solution $\phi_a$ and the BEM solution $\phi$ along the duct at 1200\,Hz with $h_1 = 0.02\tm$ is depicted, in Fig.~\ref{Fig:Anglealong}b the difference in their polar angles are shown. 
  
Compared to the behavior for the surface elements, this rise in the relative error $\frac{|\phi(\bx) - \phi_a(\bx)|}{|\phi_a(\bx)|}$ is moderate and almost not visible for smaller element widths. For finer grids the overall level of the error seems to remain constant, but with rising frequency the error level also starts to become bigger along the length of the duct.
\begin{table}[!h]
  \begin{center}
    \caption{Summary of the meshes used in Sec.~\ref{Sec:ErrorAlong} and the figures they were used for. The sidewalls $\Gamma$ are discretized using rectangular  elements of size $h_1 \times h_2$. The closures $\Gamma_0$ and $\Gamma_L$ are discretized using square elements with edge lengths $h_0$ and $h_L$. All lengths are given in m.}\label{Tab:ErrorAlong}
    \begin{tabular}{cccccc}
      Nr. & Elem. type & $h_1\um$ & $h_2\um$ & $h_0\um,h_L\um$ & \multicolumn{1}{c}{Figure} \\
      \hline
      1 & thin & 0.04 & 0.04 & 0.04 & Figs~\ref{Fig:RelErrorDiffn1}a,~\ref{Fig:RelErrorDiffn1}b,~\ref{Fig:RelErrorDiffn1}c\\
      2 & thin & 0.04 & 0.04 & 0.2 & Figs~\ref{Fig:RelErrorDiffn1}a,~\ref{Fig:RelErrorDiffn1}b,~\ref{Fig:RelErrorDiffn1}c\\
      3 & surface & 0.04 & 0.04 & 0.04 & Figs~\ref{Fig:RelErrorDiffn1}a,~\ref{Fig:RelErrorDiffn1}b,~\ref{Fig:RelErrorDiffn1}c\\
      4 & thin & 0.02 & 0.04 & 0.04 & Figs~\ref{Fig:RelErrorDiffn1}e,~\ref{Fig:RelErrorDiffn1}f,~\ref{Fig:RelErrorDiffn1}g\\
      5 & thin & 0.02 & 0.04 & 0.2 & Figs~\ref{Fig:RelErrorDiffn1}e,~\ref{Fig:RelErrorDiffn1}f,~\ref{Fig:RelErrorDiffn1}g\\
      6 & surface & 0.02 & 0.04 & 0.04 & Figs~\ref{Fig:RelErrorDiffn1}e,~\ref{Fig:RelErrorDiffn1}f,~\ref{Fig:RelErrorDiffn1}g\\
      7 & thin & 0.01 & 0.04 & 0.04 & Figs~\ref{Fig:RelErrorDiffn1}g,~\ref{Fig:RelErrorDiffn1}h,~\ref{Fig:RelErrorDiffn1}i\\
      8 & thin & 0.01 & 0.04 & 0.2 & Figs~\ref{Fig:RelErrorDiffn1}g,~\ref{Fig:RelErrorDiffn1}h,~\ref{Fig:RelErrorDiffn1}i\\
      9 & surface & 0.01 & 0.04 & 0.04 & Figs~\ref{Fig:RelErrorDiffn1}g,~\ref{Fig:RelErrorDiffn1}h,~\ref{Fig:RelErrorDiffn1}i\\
      10 & thin & 0.02 & 0.04 & 0.04 & Fig.~\ref{Fig:Anglealong}\\
    \end{tabular}
  \end{center}
\end{table}
\subsubsection{Field on the outside wall of the duct}\label{Sec:Outside}
\begin{figure}[!h]
  \begin{center}
    \includegraphics[width=0.95\textwidth]{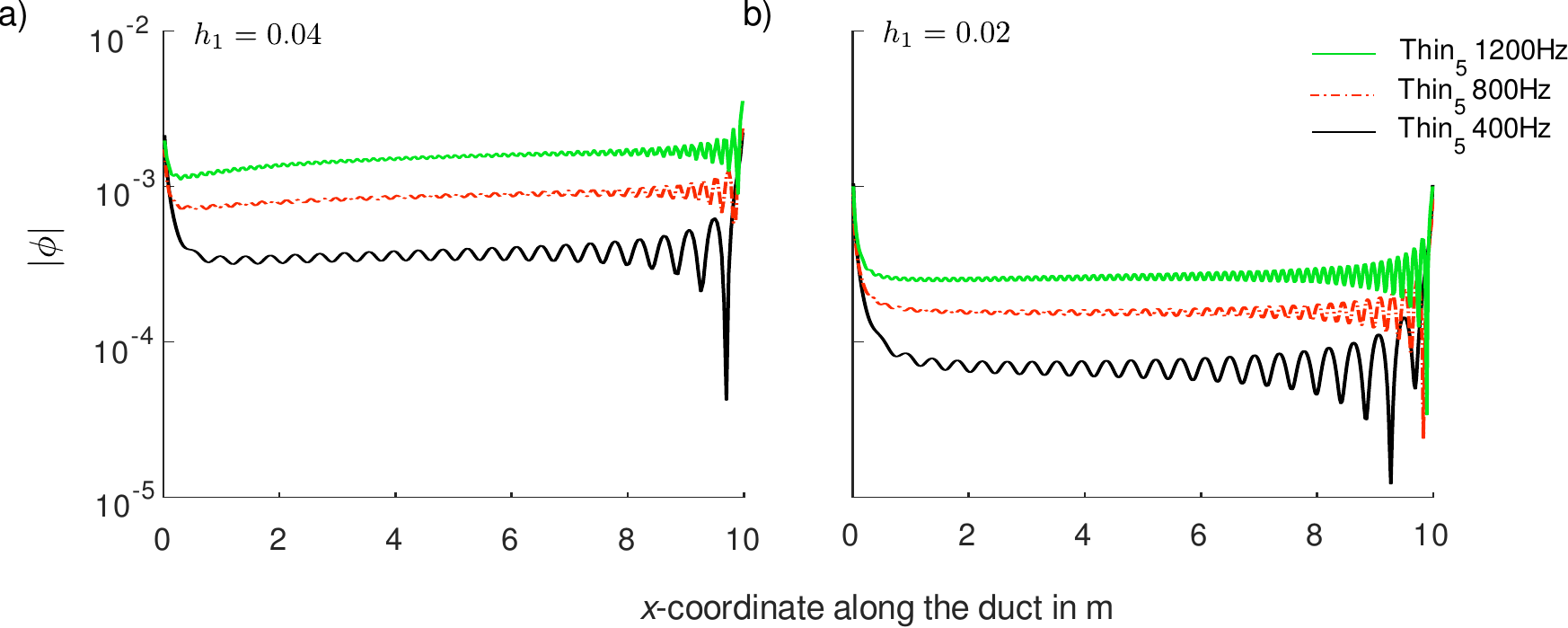}
  \end{center}
  \caption{Velocity potential $\phi^+$ on the outside along the wall at 400\,Hz, 800\,Hz and 1200\,Hz for two different widths of the boundary elements: a) $h_1 = 0.04\tm$, b) $h = 0.02\tm$. The different lines represent the velocity potential for the three frequencies $f = 400$\,Hz (black line), $f = 800$\,Hz (dash dotted line), and $f = 1200$\,Hz (light gray/green line).}\label{Fig:ErrorOutside}
\end{figure}

As already mentioned, using thin elements for this benchmark allows for a second error indicator. Obviously, the sound field outside a closed duct should be zero, thus, values of $\phi^+\neq 0$ at the outer walls of the duct are caused by numerical errors. This error is depicted in Fig.~\ref{Fig:ErrorOutside}, where the absolute value of the velocity potential $\phi^+$ on the outer part of the duct is shown for frequencies of 400\,Hz, 800\,Hz, and 1200\,Hz, and two different element widths ($h_1 = 0.04\tm$ in Fig.~\ref{Fig:ErrorOutside}a, and $h_2 = 0.02\tm$ in Fig.~\ref{Fig:ErrorOutside}b). In all cases the duct was closed at $x = 0$ and $x = L$ with elements of size $0.04\tm \times 0.04\tm$. As expected this error gets smaller with finer discretization and it is frequency dependent.

\begin{table}[!h]
  \begin{center}
        \caption{Summary of the meshes used in Sec.~\ref{Sec:Outside} and the figures they were used for. The sidewalls $\Gamma$ are discretized using rectangular thin elements of size $h_1 \times h_2$. The closures $\Gamma_0$ and $\Gamma_L$ are discretized using square thin elements with edge lengths $h_0$ and $h_L$. All lengths are given in m.}
    \begin{tabular}{ccccc}
      Nr & $h_1\um$ & $h_2\um$ & $h_0\um,h_L\um$ & Figure\\
      \hline
      1 & 0.04 & 0.04 & 0.04 & Fig.~\ref{Fig:ErrorOutside}a\\
      2 & 0.02 & 0.04 & 0.04 & Fig.~\ref{Fig:ErrorOutside}b\\
    \end{tabular}
  \end{center}
\end{table}

\section{Half Open Duct: Radiation Impedance}\label{Sec:Impedance}
The big advantage of thin boundary elements becomes clear when modeling ducts that are  open at one or both ends because no additional models for sound radiation at the duct ends are necessary. However, when using regular surface elements for thin walls, numerical problems may arise because some elements are very close to each other and the singularity of the Green's function causes numerical problems. Using the formulation for thin elements reduces some of these problems. In this example the walls of the duct are assumed to be sound hard, thus, the linear system of equations for BEM with thin elements is half the size as the one for regular surface elements.
\begin{figure}[!h]
  \begin{center}
    \includegraphics[width=0.65\textwidth]{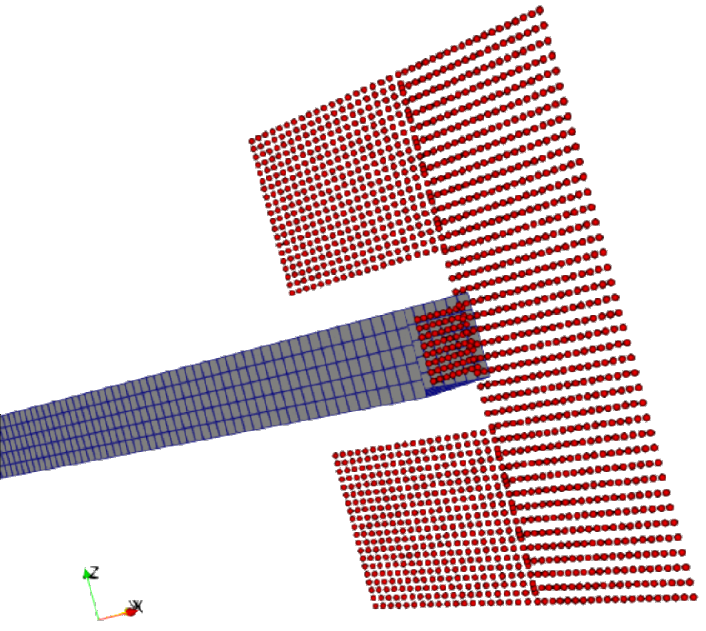}
  \end{center}
  \caption{Mesh of the half open duct with the evaluation grid at its end. The dots indicate evaluation points, where the acoustic field is calculated.}\label{Fig:OpenMesh}
\end{figure}
In the following, we look at the radiation and the radiation impedance for a 3.4\,m long square duct with width 0.2\,m using the default element size $h_1 \times h_2 = 0.04\tm \times 0.04\tm$. At $x = 0$ the duct is closed, and a velocity boundary condition $v(\bx) = 1$\,ms$^{-1}$, $\bx \in \Gamma_0$, is applied. At $x = L$ the duct is open  and $N = 64$ evaluation points $\bx_i \notin \Gamma, i = 1,\dots,N,$ are placed in the duct opening. Additionally, to illustrate the acoustic radiation from the duct, evaluation points are placed around the opening, see Fig.~\ref{Fig:OpenMesh}, where the evaluation points are depicted by small spheres. At the evaluation points, the velocity potential and the particle velocity in the $x$ direction are calculated using Green's representation theorem Eq.~(\ref{Equ:GreensRepresentation}). 
As long as the evaluation point is sufficiently far away from the boundary, the numerical evaluation of the integrals involving $G(\bx_i,\by)$ and $H(\bx_i,\by)$ does not pose a challenge. 
\begin{figure}[!h]
    \begin{center}
      \includegraphics[width=0.95\textwidth]{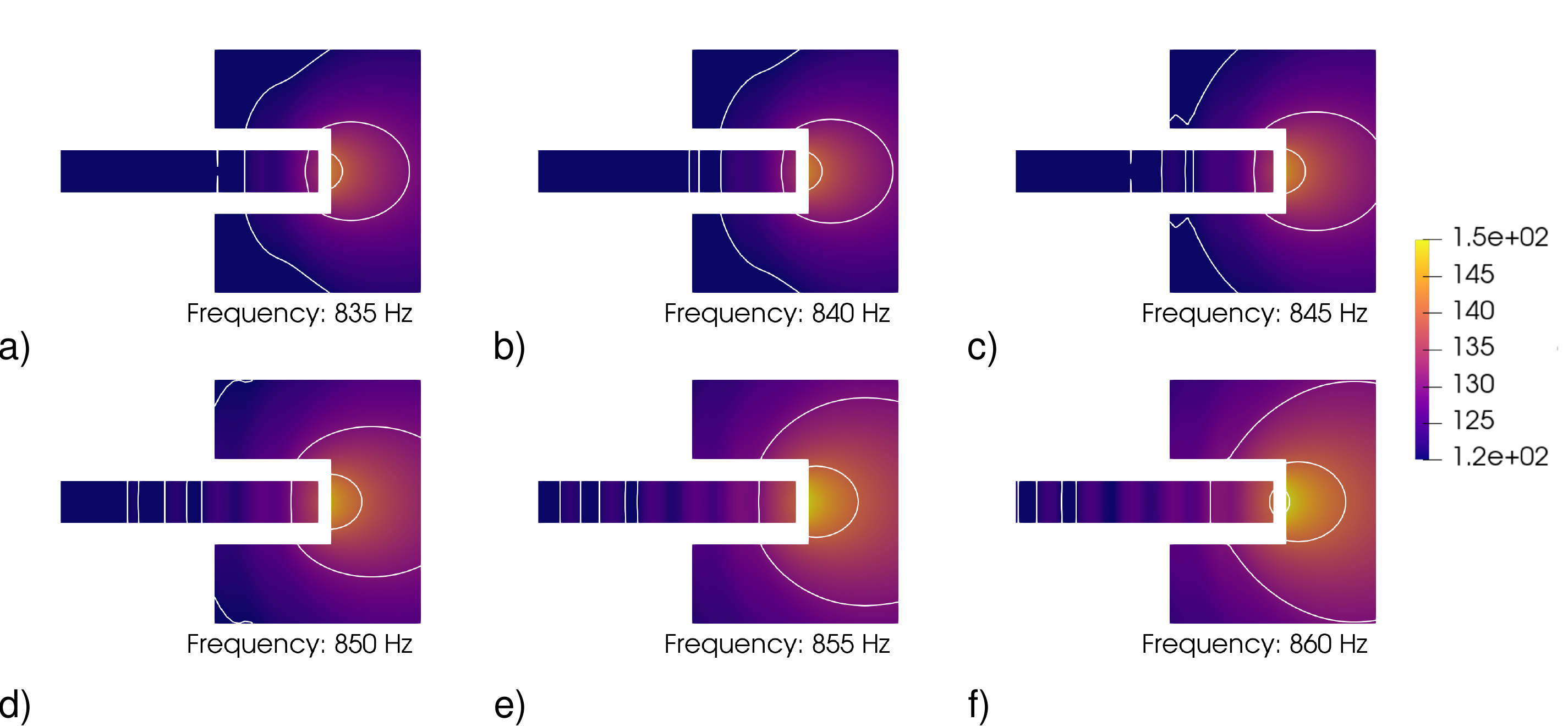}
    \end{center}
    \caption{Sound pressure at the end of the half open duct in dB between $f = 835$\,Hz and $f = 860$\,Hz. Contour lines are drawn every 10\,dB.}\label{Fig:RadField}
  \end{figure}
In Fig.~\ref{Fig:RadField} the sound pressure in dB around the duct is depicted in steps of 5\,Hz between  835\,Hz, where the mean acoustic pressure over the evaluation grid is close to a minimum, and 860\,Hz, where the mean acoustic pressure is close to a maximum. Contour lines are depicted between 120\,dB and 150\,dB in 10\,dB steps. 

\begin{figure}[!h]
  \begin{center}
    \includegraphics[width=0.75\textwidth]{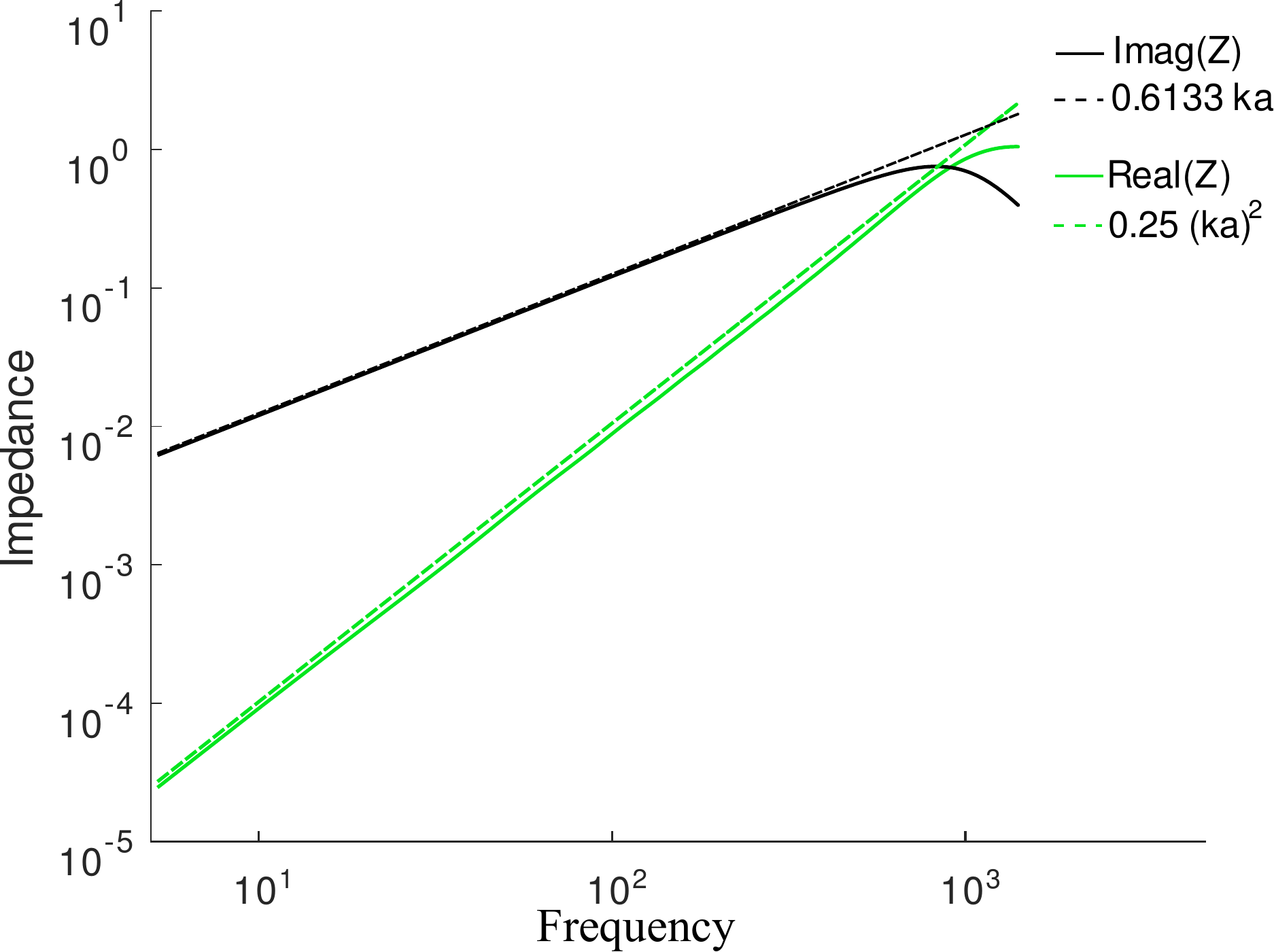}
  \end{center}
  \caption{Impedance $Z$ at the end of the duct as a function of the frequency. The dashed lines show values for an approximation formula for $Z$  used in the literature as a function of the wavenumber $k$ and the equivalent radius $a = \frac{w}{\sqrt{\pi}}$.}\label{Fig:Impedance}
\end{figure}

  In Fig.~\ref{Fig:Impedance} the relatively to air scaled impedance $Z = \frac{1}{\rho c}\frac{p}{v}$ at the open end is plotted in steps of 5\,Hz from 5\,Hz to 1415\,Hz.   The impedance at $x = L$ is calculated as the mean value over all evaluation points at $x = L$:
  $$
  Z(L) = \frac{1}N \sum\limits_{i = 1}^N \frac{p(\bx_i)}{v(\bx_i)}.
  $$
 The calculated impedance is compared with the values of the  commonly used approximation formula  for the radiation of an unflanged (circular) duct found in literature $\text{Re}(Z) \approx 0.25 (ka)^2$ and $\text{Im}(Z) \approx 0.6 ka$, cf. \cite{FleRos91,LevSch48,Creasy16}. In this formula $k$ denotes the wavenumber and  $a = \frac{w}{\sqrt{\pi}}$ is given by the equivalent radius of a circular tube with the same cross section as the benchmark duct. 
\section{Open Duct: Resonance Frequencies}\label{Sec:Boomwhacker}
In this example the sound radiation from a musical tube (boomwhacker, Klangr\"ohre, sound tube), that is a simple percussive plastic toy, will be investigated. The tube in question is approximately $63.3$\,cm long, has a radius of about 1.6\,cm,  and is ``tuned'' to the musical note C$_4 \approx 261.63$\,Hz. The length and the radius of the tubes may vary from manufacturer to manufacturer, for example, the tubes used in \cite{Ruiz14} had slightly different measurements.

\begin{figure}[!h]
  \begin{center}
    \includegraphics[width=0.45\textwidth]{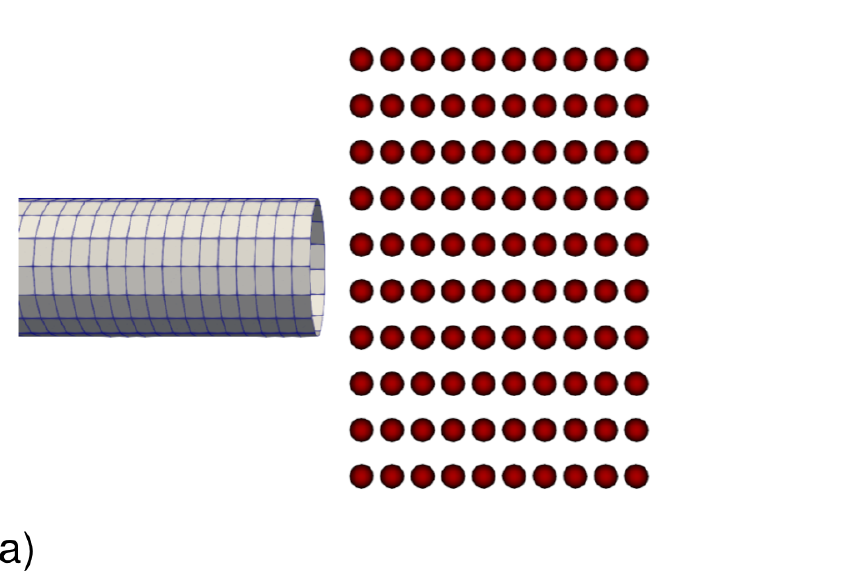}
    \includegraphics[width=0.5\textwidth]{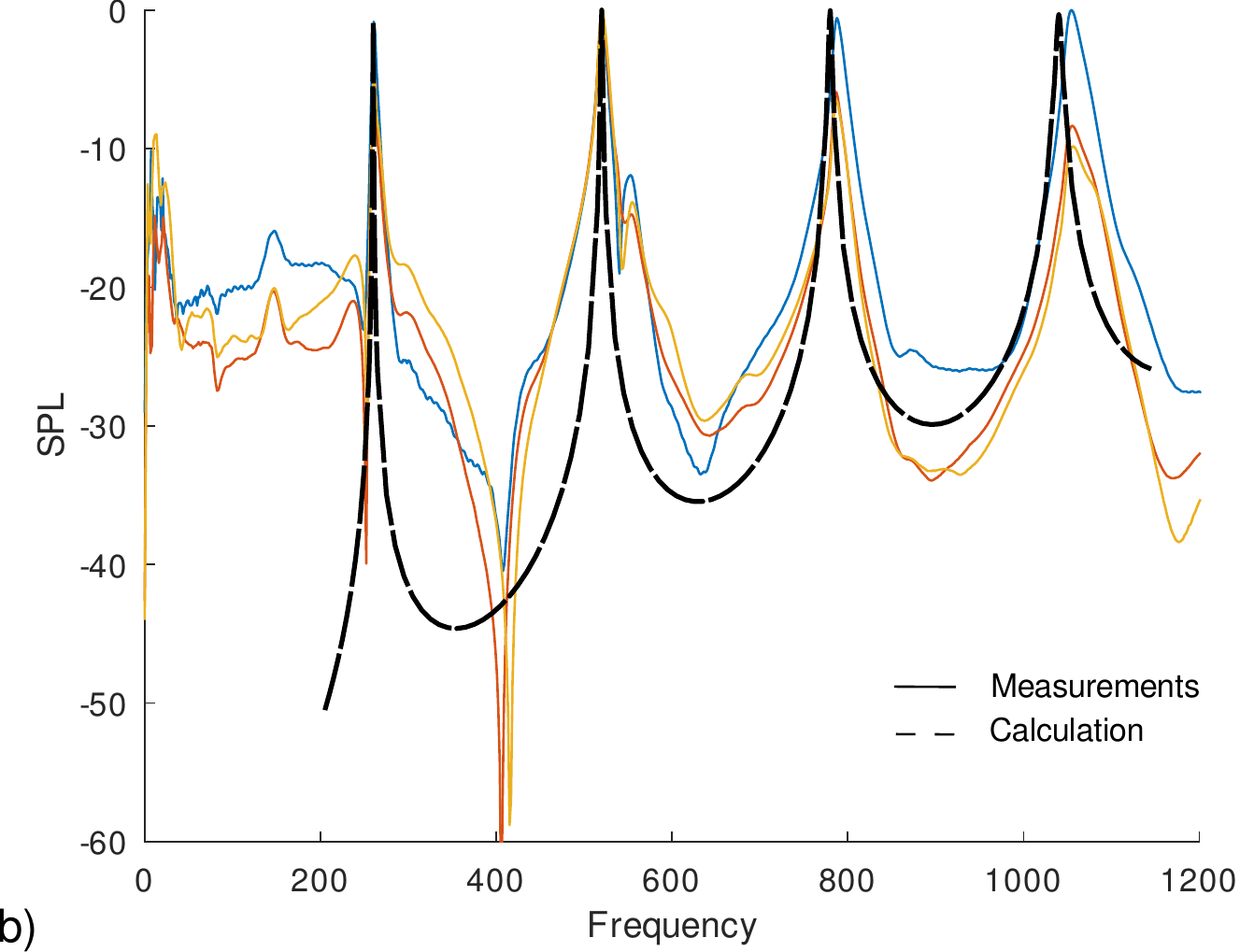}
    \caption{a) Evaluation grid at one end of the circular duct. The duct is discretized using elements of size $h_1 \times h_2 \approx 3.981\,\text{mm} \times 6.243$\,mm. b) Mean calculated pressure in dB over the evaluation grid (dashed line) and SPL for 3 measurements (continuous lines). For the calculation the mesh on the left side was used.}\label{Fig:Boomwhackerresults}
  \end{center}
\end{figure}

For the numerical calculations BEM with thin elements is used, in the different experiments (BEM30, BEM80, BEM120, and BEM160) the width of each element is given by $h_1\approx 21.82$\,mm,  $8.015$\,mm, $5.319$\,mm, and $3.981$\,mm, respectively. For all BEM calculations the height of each element was approximately $h_2 \approx 6.243$\,mm, see  Fig.~\ref{Fig:Boomwhackerresults}a for the mesh with $h_1 \times h_2 \approx 3.981\,\text{mm} \times 6.243$\,mm. Close to the open end at $x = L$ a grid of $10\times 10$ evaluation points has been placed. At a frequency of $f = 1400$\,Hz, the acoustic wavelength $\lambda = \frac{c}{f}$ is about $\lambda \approx 24.3$\,cm, thus, the general rule of thumb of about 6-10 elements per wavelength is fulfilled for all tested $h_1$ settings.
In Fig.~\ref{Fig:Boomwhackerresults}b the mean acoustic pressure over the evaluation grid at one end of the sound tube is depicted as function of frequency  for three different measurements and the BEM calculations for $h_1 \approx 3.981\,$mm.

In the calculations the minimum frequency step around the resonance frequencies is set to 1\,Hz, which means that the peaks in the BEM curve in Fig.~\ref{Fig:Boomwhackerresults} will not be exactly at the resonance frequencies of the system, and thus, have  a finite value. In the simulations, hitting the tube to produce a sound was modeled by a velocity boundary condition at an element about 3\,cm away from one end. For the measurements the tube was hit at about the same position. The produced sound was recorded with a sampling rate of $f_s = 48000$\,Hz, the length of the audible sound was about 0.6\,s, thus, the frequency resolution in the measurements is bigger than 1\,Hz. To enhance this resolution the recordings have been post-processed with the large time-frequency analysis toolbox \texttt{LTFAT} (https://ltfat.org) with approaches also used to, e.g., reassign spectrograms \cite{PruHol16}.

Especially at the first two resonances the measurements and the calculation BEM160 agree quite well, for the two higher ``harmonics'' there is a difference of a few Hz. It has to be mentioned that with higher frequencies the difference in frequency between two musical notes becomes larger, so the 4th resonance of the measurement (about 1054\,Hz) and the calculation (1040\,Hz) would register both as C$_6$ (1046.5\,Hz). The difference between measurement and calculation may be explained by the fact, that for the calculation effects at the duct wall like viscosity or heat losses have been neglected. 

The calculated peak frequencies for different discretizations are compared in Tab.~\ref{Tab:Resonances}, where only the rounded mean values of the peak values over all 3 measurements has been used. The results in the table show that the usual rule of about 6 - 10 elements per wavelength (BEM30) is not enough to get resonance frequencies accurately enough and that finer grids along the length of the duct are needed.

\begin{table}[!h]
  \begin{center}
        \caption{Resonance peaks ($F_1, F_2, F_3$ and $F_4$) of the BEM calculations and the rounded mean resonance frequency of the  measurements.}\label{Tab:Resonances}
        
    \begin{tabular}{lcccc}
      & $F_1$\,[Hz] & $F_2$\,[Hz] & $F_3$\,[Hz] & $F_4$\,[Hz]\\
      \hline
      Measurements & 261 & 520 & 788 & 1056 \\
      BEM30 ($h_1  \approx 2.18$\,cm) & 253 & 506 & 760 & 1013\\
      BEM80 ($h_1 \approx 8.015$\,mm) & 259 & 518 & 777 & 1036\\
      BEM120 ($h_1 \approx 5.319$\,mm) & 260 & 519 & 779 & 1039\\
      BEM160 ($h_1 \approx 3.981$\,mm) & 260 & 520 & 780 & 1040\\
    \end{tabular}
 \end{center}
\end{table}

\section{Conclusions}\label{Sec:Summary}
In this work, a 3D model for open and closed ducts has been presented that is based on the BEM with thin elements. These type of elements have the advantage that very thin structures can be modeled efficiently, and that compared to other numerical methods like FEM the acoustic field outside an open duct can be calculated very easily and efficiently. 

The numerical properties of the thin elements were investigated using the benchmark problem of a closed square duct \cite{Hornikxetal15}, for which the analytic solution is given by one single plane wave traveling along the duct. For this example, a slightly different behavior between BEM with regular surface elements and BEM with thin elements has been found. In general, there is greater variation of the amplitude of the velocity potential when using thin elements as compared to surface elements. One explanation for this variation may be, that the impedance boundary condition at one end of the duct is not resolved properly. In the benchmark problem the impedance $Z(L) = \rho c$ was specifically chosen because for this impedance the solution only consists of one plane wave traveling along the duct. If this specific value for the impedance is now perturbed, the analytic solution contains a second plane wave, which causes a variation in the solution.

However, compared to regular surface elements, a rise in the difference between analytic and calculated solution along the length of the duct is less prominent for thin elements, see for example Fig.~\ref{Fig:RelErrorDiffn1}. Also, a damping of the amplitude along the duct that can be observed for surface elements does not occur for thin elements. 

For both, surface and thin elements, a ``local'' behavior of the solutions can be observed, see for example Figs.~\ref{Fig:DiffDisc} and \ref{Fig:Midfacesol}. Different mesh sizes in different parts of the duct influence the solution locally in these parts. Thus, in order to achieve accurate solutions over the whole length of the duct, it is recommended to choose a \emph{uniform} element size along the whole length of the duct. 

The main advantage of thin elements was demonstrated in Sec.~\ref{Sec:Impedance} using the example of a duct that is open at the end. In this case, the advantage of the BEM that only the surface of the duct needs to be discretized to calculate the field around the object becomes apparent. As sound hard boundary conditions were used, the linear system of equations only consists of $N$ unknowns, where $N$ is the number of elements, making the formulation for thin elements efficient. When comparing the impedance at the open duct end with values found in literature, a good agreement can be found for the investigated range of frequencies.

In a third example, the spectrum of a sound tube, that is open at both ends and that is tuned to the note 'C',  was calculated and compared with measured spectra of the same tube. If the discretization is fine enough, the frequencies of the first first two resonances agree quite well, for higher harmonics there is a difference in measured and calculated peaks in the spectrum. In all calculations the resonance peaks have a strict harmonic structure, i.e., all resonances are an integer multiple of a fundamental frequency. This is not the case in the measurements. Also the peaks around the resonances are much sharper in the calculation. 

This indicates, that simply solving the Helmholtz equation is not enough especially in higher frequency bands. Realistic models will need to include losses along the duct, which can be for example modeled by a complex valued wavenumber, by using impedance boundary conditions on the sidewalls of the duct, or by explicitly modeling viscosity effects close to the duct walls. Also, the numerical experiments in this paper suggest that the common 8 to 10 elements per wavelength rule may not yield accurate results when modeling ducts.

In summary, although wave propagation inside a duct is a standard example for the 1D wave equation, there is still need to investigate realistic but efficient models for simulating the acoustic of ducts, which will be left to future work.
\section*{Declaration of Interest}
This research did not receive any specific grant from funding agencies in the public, commercial, or not-for-profit sectors.

\bibliographystyle{elsarticle-num} 
\bibliography{Kreuzer21}{}

\end{document}